\documentclass{bmcart}

\usepackage{amsthm, amsmath, amssymb}
\usepackage[per-mode=reciprocal]{siunitx}
\usepackage{booktabs}
\usepackage[section]{placeins}
\usepackage[noline, linesnumbered, noend, ruled]{algorithm2e}
\usepackage{subcaption}
\usepackage{graphicx}
\usepackage{nicefrac}
\usepackage{braket}
\usepackage{cite}

\usepackage{rotating}
\newcommand{\sw}[1]{\begin{sideways}#1\end{sideways}}

\RequirePackage{hyperref}
\usepackage{hyperref}
\usepackage[utf8]{inputenc} 

\startlocaldefs

\newcommand{\dx}[1]{\left.\text{d}#1 \right.}

\newcommand{\R}{\mathbb{R}}
\newcommand{\normal}{n}

\newcommand{\subcool}{\text{cool}}
\newcommand{\subrad}{\text{rad}}
\newcommand{\subamb}{\text{amb}}
\newcommand{\subblood}{\text{b}}
\newcommand{\subnative}{n}
\newcommand{\subcoag}{c}

\newcommand{\temperature}{T}
\newcommand{\radiation}{\varphi}
\newcommand{\damage}{\omega}

\newcommand{\conductivity}{\kappa}
\newcommand{\density}{\rho}
\newcommand{\heatcapacity}{C_p}
\newcommand{\perfusion}{\xi}

\newcommand{\absorption}[1]{\mu_{a#1}}
\newcommand{\scattering}[1]{\mu_{s#1}}
\newcommand{\anisotropy}[1]{g_{#1}}
\newcommand{\timehorizon}{\tau}
\newcommand{\htc}{\alpha}
\newcommand{\diffusivity}{D}

\newcommand{\coolingfactor}{\beta_q}

\newcommand{\frequencyfactor}{A}
\newcommand{\activationenergy}{E_a}
\newcommand{\gasconstant}{R}

\newenvironment{eqsystem*}
{\begin{equation*}
	\left\lbrace
	\begin{alignedat}{2}
}
{\end{alignedat}
	\right.
	\end{equation*}
}

\endlocaldefs

\begin{document}

\begin{frontmatter}

\begin{fmbox}
\dochead{Research}


\title{Mathematical Modeling of Vaporization during Laser-Induced Thermotherapy in Liver Tissue}

\author[
	addressref={itwm, tuk},
	email={sebastian.blauth@itwm.fraunhofer.de}
]{\inits{S}\fnm{Sebastian} \snm{Blauth}}
\author[
	addressref={kgu}, 
	email={frank.huebner@kgu.de}
]{\inits{M}\fnm{Frank} \snm{Hübner}}
\author[
	addressref={itwm},
	corref={itwm},
	email={christian.leithäuser@itwm.fraunhofer.de}
]{\inits{C}\fnm{Christian} \snm{Leithäuser}}
\author[
	addressref={itwm},
	email={norbert.siedow@itwm.fraunhofer.de}
]{\inits{N}\fnm{Norbert} \snm{Siedow}}
\author[
	addressref={kgu},
	email={t.vogl@em.uni-frankfurt.de}
]{\inits{T.J.}\fnm{Thomas J.} \snm{Vogl}}


\address[id=itwm]{
	\orgname{Fraunhofer Institute for Industrial Mathematics ITWM}, 
	\street{Fraunhofer Platz 1},                     %
	\postcode{67663}                                
	\city{Kaiserslautern},                              
	\cny{Germany}                                    
}
\address[id=tuk]{%
	\orgname{TU Kaiserslautern},
	\street{Gottlieb-Daimler Straße 48},
	\postcode{67663}
	\city{Kaiserslautern},
	\cny{Germany}
}
\address[id=kgu]{
	\orgname{Institute for Diagnostic and Interventional Radiology of the J.W. Goethe University Hospital},
    	\street{Theodor-Stern-Kai 7},
	\postcode{60590}
	\city{Frankfurt/Main}
	\cny{Germany}
}


\begin{artnotes}
\end{artnotes}

\end{fmbox}


\begin{abstractbox}

\begin{abstract} 
Laser-induced thermotherapy (LITT) is a minimally invasive method causing tumor destruction due to heat ablation and coagulative effects. Computer simulations can play an important role to assist physicians with the planning and monitoring of the treatment. Our recent study with ex-vivo porcine livers has shown that the vaporization of the water in the tissue must be taken into account when modeling LITT. We extend the model used for simulating LITT to account for vaporization using two different approaches. Results obtained with these new models are then compared with the measurements from the original study.

\end{abstract}


\begin{keyword}
	\kwd{LITT}
	\kwd{bio-heat equation}
	\kwd{vaporization}
	\kwd{enthalpy method}
	\kwd{simulation}
\end{keyword}

\end{abstractbox}
%

\end{frontmatter}




\section{Introduction}
\label{sec:introduction}

Thermal ablation methods briefly generate cytotoxic temperatures in tumorous tissue in order to destroy it. These minimally invasive methods are used for treating cancer, e.g., in lung, liver, or prostate, when surgical resection is either not possible or too dangerous for the patient. All of these methods utilize the fact that tumorous tissue is more susceptible to heat than healthy tissue to destroy as little healthy tissue as possible. Among the most common thermal ablation methods are LITT, radio-frequency ablation, and microwave ablation.

The principle of LITT \cite{muller1995laser} is based on the local supply of energy via an optical fiber, located in a water-cooled applicator. This applicator is placed directly into the tumorous tissue. The LITT treatment can take place under MRI control because the laser applicator is sourced by an optical fiber and does not include any metal parts. Therefore the patient is not exposed to radiation, in contrast to other treatments that can only be carried out under CT control.

For the therapy planning of LITT, accurate numerical simulations are needed to guide the practitioner in deciding when to stop the treatment. Mathematical models for this have been proposed, e.g., in \cite{Fasano, Mohammed}. The liver consists of about 80~\% water which vaporizes if the temperatures during the treatment become sufficiently large. The vaporization of this water is currently not included in these models but our study in \cite{Huebner, Leit} suggests that this effect is relevant for an accurate simulation. In this study the ex-vivo experiments with a larger power of 34 W show a good agreement between measured and simulated temperature until the temperature reaches approximately 100~\si{\celsius}. Then, the measured temperature stagnates while the simulated one rises further (cf. \cite{Huebner}, Fig. 3). We presume that this happens due to phase change of water which was not included in the model we used.

In this paper we use the measurements from \cite{Huebner} and compare two models for the vaporization. One of them is the effective specific heat (ESH) model introduced in \cite{Yang} which modifies the specific heat coefficient to account for the phase change. The other one is the enthalpy model which uses an additional state equation to model the phase transition. We compare the models to experimental data with ex-vivo porcine livers from \cite{Huebner}. When modeling the vaporization it is also necessary to consider what happens to the vapor. We do this with a very simple condensation model also proposed in \cite{Yang}, where no transport mechanism is involved for the vapor. Instead it is assumed that it condensates instantaneously in a region with lower temperature. The effects of this simplification are discussed.

This paper is structured as follows. 
Our existing mathematical model for simulating LITT including heat and radiative transfer is described in Section~\ref{sec:mathematical_model}. This model is based on the work of \cite{Fasano} and we have also used it in \cite{Huebner}. In Section~\ref{sec:vapo} we modify and extend this model in such a way that it also covers the effect of vaporization during the treatment. Therefore, we consider both the ESH model of \cite{Yang} as well as an enthalpy model for vaporization.
Afterwards, we present the details of the numerical solution of our models in Section~\ref{sec:num_sim}. Finally, the models are validated with measurement data obtained from experiments made with ex-vivo porcine liver tissue (cf. \cite{Huebner}) in Section~\ref{sec:results}.


\section{Mathematical Model}
\label{sec:mathematical_model}

\begin{figure}[h!bt]
	\centering
	\includegraphics[width=0.5\textwidth]{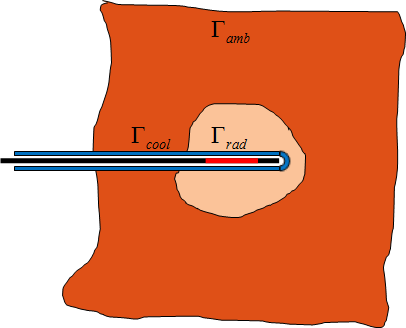}
	\centering
	\caption{Sketch of the geometry including the water-cooled applicator with radiating laser fiber.}
	\label{fig:sketch_geometry}    
\end{figure}

We denote by $\Omega \subset \R^3$ the geometry of the liver and by $\Gamma = \partial \Omega$ its boundary. The latter consists of the radiating surface of the adjacent applicator $\Gamma_\subrad$, the cooled surface of the applicator $\Gamma_\subcool$, and the ambient surface of the liver $\Gamma_\subamb$ (see Figure~\ref{fig:sketch_geometry}).
The mathematical model is described by a system of partial differential equations (PDEs) for the heat transfer inside the liver, the radiative transfer from the applicator into the liver tissue, and a model for tissue damage (cf. \cite{Mohammed, Fasano, Huebner}).  

\subsection{Heat Transfer}
The heat transfer in the liver tissue is modeled by the well-known {\it bio-heat equation} (cf. \cite{pennes})
\begin{alignat}{2}
	\label{eq:bioheat}
	\density \heatcapacity \frac{\partial \temperature}{\partial t} - \nabla\cdot\left(\conductivity \nabla \temperature\right) + \perfusion_\subblood(\temperature - \temperature_\subblood) &= Q_\subrad \quad &&\text{ in } (0, \timehorizon) \times \Omega, \\
	\temperature(0,\cdot) &= \temperature_\text{init} \quad &&\text{ in } \Omega, \nonumber
\end{alignat}
where $\temperature = \temperature(x,t)$ denotes the temperature of the tissue, depending on the position $x\in \Omega$ and the time $t\in (0,\timehorizon)$. Here, the end time of the simulation is denoted by $\timehorizon > 0$. Further, $\heatcapacity$ is the specific heat capacity, $\density$ the density of the tissue, and $\conductivity$ the thermal conductivity. The perfusion rate due to blood flow is denoted by $\perfusion_\subblood$ and the blood temperature by $\temperature_\subblood$. Note that in the current ex-vivo study the perfusion rate $\perfusion_\subblood$ is set to zero. Finally, $Q_\subrad$ is the energy source term due to the irradiation of the laser fiber and the initial tissue temperature distribution is given by $\temperature_\text{init}$.

For the heat transfer between the tissue and its surroundings, given by the ambient surface and the applicator, the following Robin type boundary conditions are used
\begin{alignat*}{2}
	\conductivity\ \partial_\normal \temperature =&\ \htc_\subcool(\temperature_\subcool - \temperature) \quad &&\text{ on } (0, \timehorizon) \times \left(\Gamma_\subrad\cup \Gamma_\subcool \right), \nonumber \\
	\conductivity\ \partial_\normal \temperature =&\ \htc_\subamb(\temperature_\subamb - \temperature) \quad &&\text{ on } (0, \timehorizon) \times \Gamma_\subamb.
\end{alignat*}
Here, $\normal$ is the outer unit normal vector on $\Gamma$. Additionally, $\htc_\subcool$ and $\htc_\subamb$ are the heat transfer coefficients for the water-cooled part of the applicator and the surroundings of the liver, respectively. The temperature of the cooling water is denoted by $\temperature_\subcool$ and $\temperature_\subamb$ is the ambient temperature. We come back to this bio-heat equation in Section~\ref{sec:vapo}, where we modify it such that it also covers the effect of vaporization of water in the tissue. The radiative source term $Q_\subrad$ is defined in the next section by \eqref{eq:source_heat}.

\subsection{Radiative Transfer}
\label{subsec:rte}
The irradiation of laser light is modeled by the {\it radiative transfer equation}
\begin{equation}
	\label{eq:rte}
	s\cdot\nabla I+\left(\absorption{} + \scattering{} \right)I=\frac{\scattering{}}{4\pi}\int\limits_{S^2}P(s\cdot s')I(s',x) \dx{s'} \quad \text{ in } S^2 \times \Omega,
\end{equation}
where the radiative intensity $I = I(s,x)$ depends on a direction $s\in S^2$ on the (unit) sphere and the position $x\in \Omega$, and $\absorption{}$ and $\scattering{}$ are the absorption and scattering coefficients, respectively. In particular, as that radiative transfer happens significantly faster than temperature transfer, we neglect the time-dependence and use this quasi-stationary model. The scattering phase function $P(s\cdot s')$ is given by the Henyey-Greenstein term which reads (cf. \cite{niemz})
\begin{equation*}
	P(s\cdot s')=\frac{1-\anisotropy{}^2}{(1+\anisotropy{}^2-2\anisotropy{}(s\cdot s'))^{3/2}}\nonumber.
\end{equation*}
Here, $\anisotropy{} \in [-1, 1]$ is the so-called anisotropy factor that describes backward scattering for $g=-1$, isotropic scattering in case $\anisotropy{}=0$ and forward scattering for $\anisotropy{}=1$.

Due to the high dimensionality of the radiative transfer equation \eqref{eq:rte}, we use the so-called $P_1$-approximation to model the radiative energy, the details of which can be found, e.g., in \cite{Modest}. Introducing the ansatz
\begin{equation*}
	I(s,x)=\phi(x)+3s\cdot q(x),
\end{equation*}
where $q(x)=\frac{1}{4\pi}\int\limits_{S^2}I(s,x) s \dx{s}$ is radiative flux vector, one obtains the much simpler three-dimensional diffusion equation
\begin{alignat}{2}
	\label{eq:p1}
	-\nabla\cdot(\diffusivity\nabla\radiation)+\mu_a \radiation =& 0 \quad &&\text{ in } \Omega,
\end{alignat}
where $\radiation = \radiation(x)$ is the radiative energy and the diffusion coefficient $D$ is given by
\begin{equation*}
	\diffusivity=\frac{1}{3(\absorption{}+(1-\anisotropy{})\scattering{})}.
\end{equation*}

To derive the boundary conditions we use Marshak's procedure as described in, e.g., \cite{Modest}. We obtain Robin type boundary conditions
\begin{equation}
	\label{eq:marshak}
	\diffusivity \frac{\partial \radiation}{\partial \normal} = \frac{q_\text{app}}{A_{\Gamma_\subrad}} \quad \text{ on } \Gamma_\subrad, \qquad \diffusivity \frac{\partial \radiation}{\partial \normal}+b\radiation=0 \quad \text{ on } \left( \Gamma_\subcool \cup \Gamma_\subamb \right),
\end{equation}
where $q_\text{app}$ is the laser power entering the tissue and $A_{\Gamma_\subrad}$ the surface area of the radiating part of the applicator. The former can be written as
\begin{equation*}
	q_\text{app}(t) = \begin{cases}
		(1-\coolingfactor) \hat{q} & \text{ if } t_\text{on} \leq t \leq t_\text{off}, \\
		0 & \text{ otherwise, }
	\end{cases}
\end{equation*}
where $\hat{q}$ is the configured laser power and the factor $(1-\coolingfactor)$ models the absorption of energy by the coolant (cf. \cite{Huebner}). Moreover, the parameter $b$ in \eqref{eq:marshak} is given as $b=0.5$ on $\Gamma_\subamb$ and $b=0$ on $\Gamma_\subcool$.
From the numerical point of view the system given by \eqref{eq:p1} and \eqref{eq:marshak} is much easier to solve than the original system given by \eqref{eq:rte}. Finally, the radiative energy is used to define the source term for the bio-heat equation in the following way
\begin{equation}
\label{eq:source_heat}
Q_\subrad(x) = \absorption{} \radiation(x).
\end{equation}

\subsection{Tissue Damage and Its Influence on Optical Parameters}
The optical parameters $\absorption{}, \scattering{}$ and $\anisotropy{}$ are very sensitive to changes of tissue's state. In particular, once the coagulation of cells starts, their optical parameters change and, as a result, the radiation cannot enter the tissue as deeply as before. Therefore, we model the damage of the tissue as in, e.g., \cite{Mohammed, Fasano} with the help of the {\it Arrhenius law}, which is given by
\begin{equation}
	\label{eq:arrhenius}
	\damage(t,x)=\int\limits_0^t \frequencyfactor \exp\left(- \frac{\activationenergy}{\gasconstant \temperature(s, x)}\right)\dx{s},
\end{equation}
with so-called frequency factor $A$, activation energy $E_a$, and universal gas constant $R$. This describes the change of optical parameters due to coagulation in the following way
\begin{align*}
	\absorption{} =&\ \absorption{\subnative} + (1-e^{-\damage}) (\absorption{\subcoag} - \absorption{\subnative}), \\
	\scattering{} =&\ \scattering{\subnative} + (1-e^{-\damage}) (\scattering{\subcoag} - \scattering{\subnative}), \\
	\anisotropy{} =&\ \anisotropy{\subnative} + (1-e^{-\damage}) (\anisotropy{\subcoag} - \anisotropy{\subnative}),
\end{align*}
where the subscripts $\subnative$ and $\subcoag$ indicate properties of native and coagulated tissue, respectively (cf. \cite{Fasano}).

\section{Mathematical Modeling of Vaporization}
\label{sec:vapo}
Vaporization of water inside organic materials plays an important role in many different fields, e.g., in medicine or the food industry. To model the temperature distribution in such materials correctly, it is important to take the vaporization into account as a significant amount of energy is necessary for the phase transition from water to vapor. The basic principle is the following (see, e.g., \cite{Demtrder2018}). If energy in the form of heat is added to water (under constant pressure), the water's temperature increases as long as it is below the vaporization temperature, i.e., below \num{100}~\si{\celsius}. However, as soon as the water reaches this temperature, the temperature does not increase further, although heat is still added to the water. At this point, the water starts to boil and eventually vaporizes after a sufficient amount of energy was added to it. Finally, the temperature of the emerging water vapor increases again after all water has been vaporized. This happens due to the fact that the energy added to the water at its boiling point is used to change its phase and not to increase its temperature, until all water is vaporized.

In the following, we discuss two vaporization models. First, we take a look at the effective specific heat (ESH) model introduced in \cite{Yang} which uses a varying specific heat capacity to model the phase change. In this model the phase transition is spread over a reasonably small interval around \SI{100}{\celsius}. This simplification makes it possible to model the phase transition using a single PDE. Second, we propose an enthalpy model with an additional state equation for the enthalpy. For this model, the transition happens at a single temperature, namely at \SI{100}{\celsius}.

\subsection{The Effective Specific Heat (ESH) Model}
\label{subsec:eff}
The ESH model introduced in \cite{Yang} considers the following modified bio-heat equation
\begin{equation}
	\label{eq:heat_yang}
	\density \heatcapacity \frac{\partial \temperature}{\partial t} - \nabla\cdot\left(\conductivity \nabla \temperature\right) + \perfusion_\subblood(\temperature - \temperature_\subblood) = Q_\subrad - Q_\text{vap} + Q_\text{cond} \quad \text{ in } (0, \timehorizon) \times \Omega, \\
\end{equation}
with the same initial and boundary conditions as \eqref{eq:bioheat}. Here, $Q_\text{vap}$ is a source term that models the vaporization of water and $Q_\text{cond}$ is the source term for the condensation (see Section~\ref{sec:condensation esh}). In \cite{Yang} this has the following form
\begin{equation}
	\label{eq:source_yang}
	Q_{vap} = - \lambda \frac{d W}{d t}
\end{equation}
where $\lambda$ denotes the latent heat of water and $W$ is the tissue water density. Using the chain rule we see that
\begin{equation*}
	\frac{d W}{d t} = \frac{\partial W}{\partial \temperature} \frac{\partial \temperature}{\partial t}.
\end{equation*}
Substituting this into \eqref{eq:source_yang} and \eqref{eq:heat_yang} gives the following modified heat equation
\begin{equation*}
	\density \heatcapacity' \frac{\partial \temperature}{\partial t} - \nabla\cdot\left(\conductivity \nabla \temperature\right) + \perfusion_\subblood(\temperature - \temperature_\subblood) = Q_\subrad + Q_\text{cond} \quad \text{ in } (0, \timehorizon) \times \Omega, \\
\end{equation*}
where the effective specific heat capacity $\heatcapacity'$ is given by
\begin{equation*}
	\heatcapacity' = \heatcapacity - \frac{\lambda}{\density} \frac{\partial W}{\partial \temperature}.
\end{equation*}
Since $\frac{\partial W}{\partial \temperature} < 0$ for vaporization (the water content decreases with temperature), we have that $\heatcapacity' \geq \heatcapacity$. 

\begin{figure}[!t]
	\centering
	\begin{subfigure}{0.475\textwidth}
		\includegraphics[width=\textwidth]{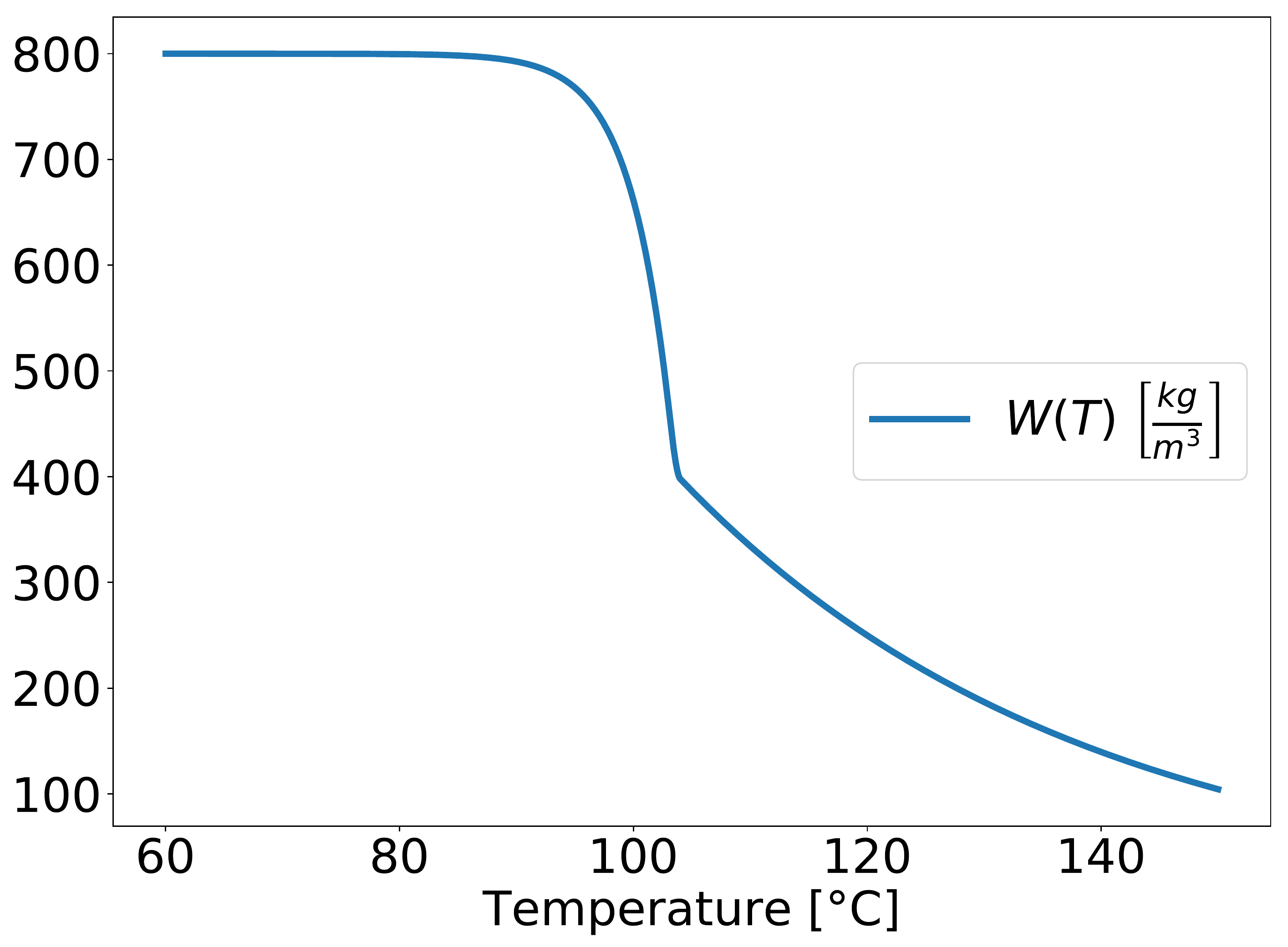}
	\end{subfigure}
	\hfil
	\begin{subfigure}{0.475\textwidth}
		\includegraphics[width=\textwidth]{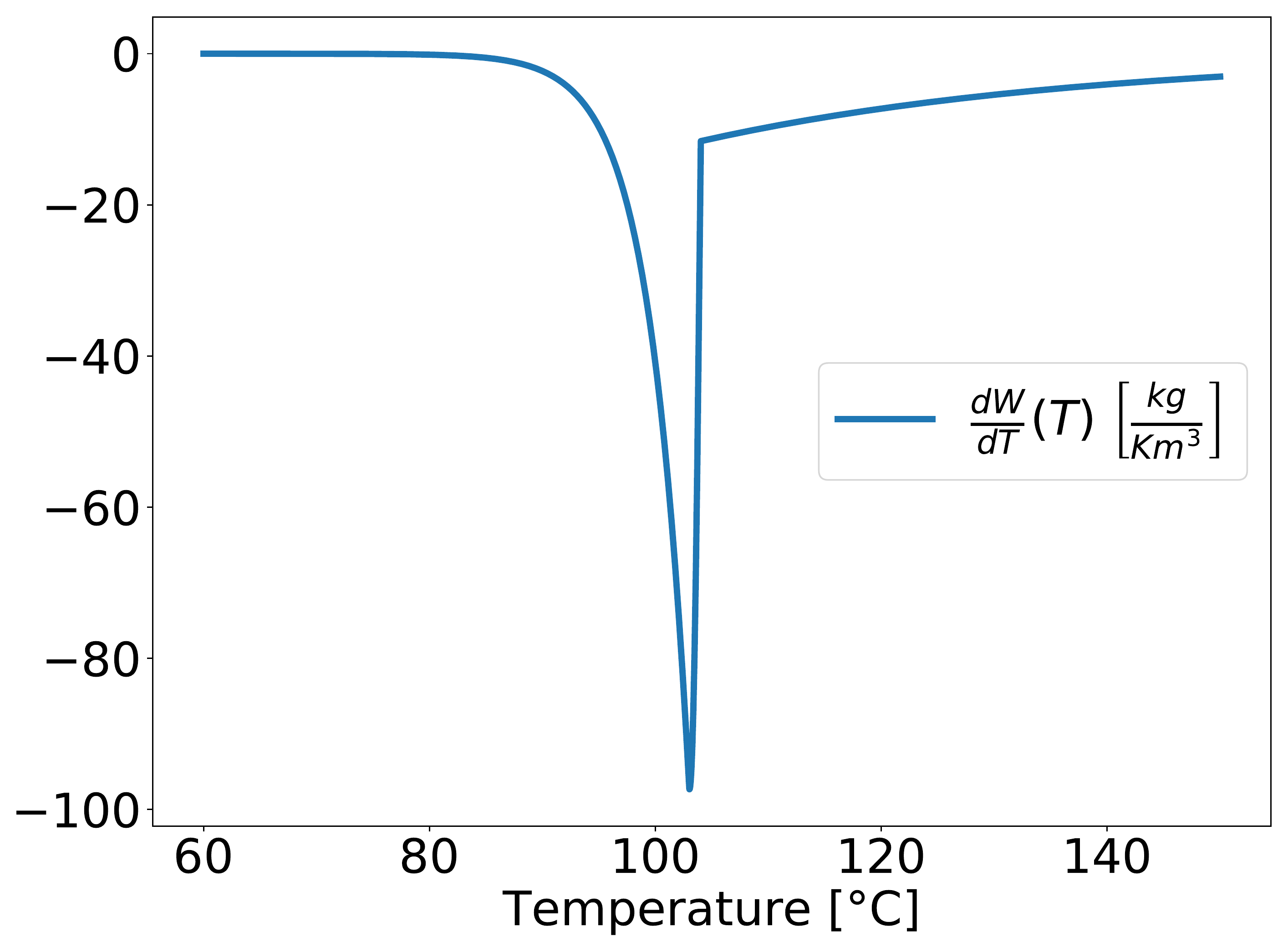}
	\end{subfigure}
	\caption{Function $W(T)$ and derivative $\frac{dW}{dT}(T)$ of tissue water density from \cite{Yang}.}
	\label{fig:2}    
\end{figure}

Based on experiments that measured water content of bovine liver as a function of temperature in \cite{Yang1} the following function is used to describe the tissue water density (cf. \cite{Yang, Yang1})
\begin{equation*}
	W(T)= 800 \cdot
	\begin{cases}
		\left( 1- e^{\frac{T-106}{3.42}} \right) & \text{ if } T\leq 103 \si{\celsius}, \\
		S(T) & \text{ if } 103 \si{\celsius} < T \leq 104 \si{\celsius},\\
		e^{-\frac{T-80}{34.37}} & \text{ if } 104 \si{\celsius} \le T,
	\end{cases}
	\label{w}
\end{equation*}
where $S(T)$ is the cubic $C^1$ spline that interpolates between the two exponential functions, (approximately) given by
\begin{equation*}
	S(T) = \num{3.712982e-2}~T^3 - \num{11.47524}~T^2 + \num{1.182046e3}~T - \num{4.058214e+04}.
\end{equation*}
The function $W$ and its derivative are depicted in Figure~\ref{fig:2}. In particular, we get that the effective specific heat is very large in an area around \SI{100}{\celsius}. Therefore, it holds that
\begin{equation*}
	\label{approx}
	\frac{\partial T}{\partial t} \ll 1 \quad \text{ for } \temperature \text{ around } \SI{100}{\celsius},
\end{equation*}
which models the vaporization of the tissue water.



\subsection{Simple Condensation Model for ESH Model}
\label{sec:condensation esh}
In \cite{Yang} it is discussed that, in addition to the vaporization of water, one also needs to consider the effect of condensation of the water vapor in order to obtain an accurate model. There, it was assumed that the water vapor diffuses into a region of lesser temperatures where it condensates and releases its latent heat obtained through the vaporization. The authors of \cite{Yang} describe their model for this in the following way. They say that they first calculate the total amount of water that was vaporized in the last time step. From this, the amount of latent heat generated is computed. Finally, this is added uniformly to the tissue region whose temperature is between \num{60}~\si{\celsius} and \num{80}~\si{\celsius}. We have implemented this simple condensation model in the following way. We compute the total amount of latent heat which is currently consumed through the vaporization of water by
\begin{equation*}
\bar{Q}_\text{vap} = \int_\Omega Q_\text{vap} \dx{x},
\end{equation*}
where [$\bar{Q}_\text{vap}$] = \si{\watt}. Additionally, we define the condensation region as
\begin{equation*}
	\Omega_\text{cond} := \Set{x\in \Omega | \temperature^- \leq \temperature \leq \temperature^+ },
\end{equation*}
for given temperature boundaries $\temperature^- < \temperature^+ < \SI{100}{\celsius}$. Uniformly distributing $\bar{Q}_\text{vap}$ over the condensation region then yields the condensation source term
\begin{equation*}
	Q_\text{cond}(x) = \begin{cases}
	\frac{\bar{Q}_\text{vap}}{\text{vol}(\Omega_\text{cond})} \quad &\text{ if } x\in \Omega_\text{cond} \text{ and } \text{vol}(\Omega_\text{cond}) >0, \\
0 \quad &\text{ otherwise}.
	\end{cases}
\end{equation*}
In particular, this implies that our model is energy conserving. This is of course a very rough condensation model because there is no real transport mechanism for the vapor involved at all. Any vapor will instantaneously condensate in another region with lower temperature. This simple model shows promising results but there is also room for improvement as discussed in Section~\ref{sec:discussion condensation model}.

\subsection{Enthalpy Model}
\label{subsec:itwm}

In the this section, we present the details of the second model for vaporization, which is based on an enthalpy formulation. It consists of two coupled equations, one for the temperature of the tissue and one for its enthalpy. For the temperature, we have the following, modified bio-heat equation
\begin{equation}
	\label{eq:bio_vap}
	\density \heatcapacity \frac{\partial \temperature}{\partial t} = \begin{cases}
	\nabla\cdot (\conductivity \nabla \temperature) + \perfusion(\temperature_\subblood - \temperature) + Q_\subrad + Q_\text{cond} & \quad \parbox{3cm}{if $\temperature < \num{100}~\si{\celsius}$ or \\ $\temperature \geq \num{100}~\si{\celsius}$ \\ and $H = \rho \lambda_\text{vap}$,}  \\
	\mbox{} \\
	0 \quad & \parbox{3.2cm}{ if $\temperature = \num{100}~\si{\celsius}$ \\ and $0 \leq H < \density \lambda_\text{vap}$,} 
	\end{cases}
\end{equation}
where $\lambda_\text{vap} = \num{0.8} \lambda$ is the proportion of the enthalpy of vaporization corresponding to the tissue's water content of \SI{80}{\percent}.
Further, the (volumetric) enthalpy of the water $H$, [$H$] = \si{\joule \per \cubic\meter}, is modeled by the following ODE
\begin{equation}
	\label{eq:enthalpy}
	\frac{\partial H}{\partial t} = \begin{cases}
	0 & \quad \parbox{5cm}{if $\temperature < \num{100}~\si{\celsius}$ or \\ $\temperature \geq \num{100}~\si{\celsius}$ and $H = \rho\lambda_\text{vap}$,} \\
	\mbox{} \\
	\nabla \cdot (\conductivity \nabla \temperature) + \perfusion(\temperature_\subblood - \temperature) + Q_\subrad & \quad \parbox{5cm}{if $\temperature = \num{100}~\si{\celsius}$ \\ and $0 \leq H < \rho\lambda_\text{vap}$.}
	\end{cases}
\end{equation}
Equation \eqref{eq:bio_vap} has the same initial and boundary conditions as \eqref{eq:bioheat}, and the initial condition of the enthalpy is given by $H = 0$ in $\Omega$, i.e., no vaporization had happened before the treatment. The term $Q_\text{cond}$ describes a heat source due to the condensation of water vapor in regions with temperatures below \num{100}~\si{\celsius}, similar to the one of the ESH model (cf. Section~\ref{sec:condensation esh}). Observe that the modified bio-heat equation \eqref{eq:bio_vap} coincides with the classical bio-heat equation \eqref{eq:bioheat} and we also have $H=0$, i.e., no vaporization is happening, as long as we have that $\temperature < \num{100}~\si{\celsius}$ everywhere. This changes as soon as $\temperature = \num{100}~\si{\celsius}$ at some point $\bar{x}\in \Omega$. Then, we see that the bio-heat equation \eqref{eq:bio_vap} gives $\frac{\partial \temperature}{\partial t}(\bar{x}) = 0$ and, therefore, $\temperature(\bar{x}) = \num{100}~\si{\celsius}$ in case $0 \leq H(\bar{x}) < \density\lambda_\text{vap}$, i.e., the temperature at a point does not change until the enthalpy exceeds the enthalpy of vaporization $\rho \lambda_\text{vap}$.
In the meantime, the energy that would usually lead to an increase in temperature now only increases the enthalpy, which models the phase change of the water in the tissue. Finally, as soon as the enthalpy reaches the enthalpy of vaporization, all water is vaporized and the bio-heat equation is valid again.

\subsection{Simple Condensation Model for Enthalpy Model}
\label{sec:condensation enthalpy}
Similar to Section~\ref{sec:condensation esh} the simple condensation model suggested in \cite{Yang} is used. In contrast to the ESH model, the total amount of latent heat can be computed from the change of enthalpy in the following way
\begin{equation}
	\label{eq:Qbarvap enthalpy}
	\bar{Q}_\text{vap} = \int_\Omega \frac{\partial H}{\partial t} \dx{x}.
\end{equation}
Again, the condensation region is defined by
\begin{equation*}
	\Omega_\text{cond} = \Set{x\in \Omega | \temperature^- \leq \temperature \leq \temperature^+ },
\end{equation*}
and the condensation source term is
\begin{equation*}
	Q_\text{cond}(x) = \begin{cases}
	\frac{\bar{Q}_\text{vap}}{\text{vol}(\Omega_\text{cond})} \quad &\text{ if } x\in \Omega_\text{cond} \text{ and } \text{vol}(\Omega_\text{cond}) >0, \\
0 \quad &\text{ otherwise},
	\end{cases}
\end{equation*}
where $\text{vol}(\Omega_\text{cond})$ denotes the volume of $\Omega_\text{cond}$. With this, we get that the temperature increase due to condensation corresponds to the energy used to change the phase of the water, uniformly distributed over $\Omega_\text{cond}$. Finally, note that the numerical discretization of this model is described in Section~\ref{sec:num_disc}.

\section{Numerical Methods}
\label{sec:num_sim}

\begin{table}[!b]
	\centering
	\begin{tabular}{|l|r|}
		\hline
		Parameter & Value \\
		\hline 
		&\\
		Optical (native): & \\
		Absorption coefficient $\absorption{\subnative}$ [\si{\per \meter}] & \num{50} \\
		Scattering coefficient $\scattering{\subnative}$ [\si{\per \meter}]  & \num{8000} \\
		Anisotropy factor $\anisotropy{\subnative}$                              & \num{0.97}\\
		\hline 
		&\\
		Optical (coagulated): & \\
		Absorption coefficient $\absorption{\subcoag}$ [\si{\per \meter}]  & \num{60} \\
		Scattering coefficient $\scattering{\subcoag}$ [\si{\per \meter}] & \num{30000} \\
		Anisotropy factor $\anisotropy{\subcoag}$                                & \num{0.95}\\
		\hline 
		&\\
		Thermal: &\\
		Heat conductivity $\conductivity$ [\si{\watt \per \meter \per \kelvin}]    &   \num{0.518} \\
		Heat capacity $\heatcapacity$ [\si{\joule \per \kilogram \per \kelvin}]        &  \num{3640} \\
		Tissue density $\density$ [\si{\kilogram \per \cubic \meter}]            & \num{1137} \\
		Heat transfer coefficient $\htc_\subcool$ [\si{\watt \per \square \meter \per \kelvin}] & \num{250} \\
		Heat transfer coefficient $\htc_\subamb$ [\si{\watt \per \square \meter \per \kelvin}] & \num{44} \\
		\hline
		&\\
		Damage: &\\
		Damage rate constant $\frequencyfactor$ [\si{\per \second}]      & \num{3.1e98} \\
		Damage activation energy $\activationenergy$ [\si{\joule \per \mole \per \kelvin}] & \num{6.3e5} \\
		Universal gas constant $\gasconstant$ [\si{\joule \per \mole \per \kelvin}] & \num{8.31}\\
		\hline
		&\\
		Vaporization: &\\
		Latent heat of water $\lambda$ [\si{\joule \per \kilogram}]    & \num{2257e3} \\
		\hline
	\end{tabular}
	\caption{Physical parameters for LITT in ex-vivo porcine liver tissue.}
	\label{table:parameters}
\end{table}

In this section, we detail the numerical methods used to discretize and solve the governing equations.

\subsection{Numerical Solution of the Governing PDEs}
\label{subsec:numerics}
The mathematical model for radiative heat transfer and the models for vaporization described above were used to simulate the behavior of ex-vivo porcine liver tissue during LITT. 
The computational geometry was generated using Open Cascade (Open Cascade SAS, Guyancourt, France) and the mesh was created with the help of GMSH, version 2.11.0 (cf. \cite{gmsh}). The governing equations were solved with the finite element method in Python, version 2.7, using the package FEniCS, version 2017.2 (cf. \cite{fenics, fenics_book}). For the numerical solution of the PDEs, we first \mbox{(semi-)discretize} the bio-heat equation in time using the implicit Euler method. Then, we use piecewise linear Lagrange elements for the spatial discretization of the temperature and radiative energy. The resulting sequence of linear systems was then solved with the help of PETSc (cf. \cite{petsc-user-ref}), where we used the conjugate gradient method with a relative tolerance of \num{1e-10}. Afterwards, the damage function is computed using a right-hand Riemann sum to discretize the time integral of \eqref{eq:arrhenius}.

\subsection{Discretization of the Enthalpy Model}
\label{sec:num_disc}

In the following we describe our discretization of the enthalpy model. In particular, to compute the temperature distribution from time $t$ to $t+\Delta t$ we proceed as follows. We first solve \eqref{eq:p1}  to obtain the radiative energy at $t+ \Delta t$. With this, we compute the temperature distribution at $t + \Delta t$ from \eqref{eq:bioheat}. Subsequently, we iterate over the nodes of the finite element mesh and check, whether the temperature exceeds \num{100}~\si{\celsius}. At these nodes, the temperature is set to \num{100}~\si{\celsius} and from the excess temperature we compute the corresponding increase in enthalpy. If the enthalpy surpasses the limit of $\density \lambda_\text{vap}$, we return this surplus in the form of heat to the corresponding nodes. After doing so, we integrate the (local) changes in enthalpy over $\Omega$ to compute the total change of enthalpy $\Delta H$. Therefore, we can now compute the source term $\bar{Q}_\text{vap}$ of \eqref{eq:Qbarvap enthalpy} as follows
\begin{equation*}
	\bar{Q}_\text{vap} = \frac{\Delta H}{\Delta t},
\end{equation*}
which is then used as the source term for the next time step, simulating the release of enthalpy by the condensation of the water vapor. Then, the new tissue damage is computed from \eqref{eq:arrhenius} and the procedure is continued until we reach the end time $\timehorizon$.

\section{Results and Discussion}
\label{sec:results}

We use the experiments from the study of \cite{Huebner} to test the vaporization models. In this study LITT was applied to ex-vivo porcine livers and the resulting temperature was measured with a probe. The experiment was repeated nine times with different laser powers and different flow rates for the applicator cooling system. For the study in \cite{Huebner}, the authors used the mathematical model introduced in Section~\ref{sec:mathematical_model} which was derived from the one presented in \cite{Fasano}. However, the model did not take into account the vaporization of water in the tissue. While the general agreement between experiment and simulation was good, there were notably two outliers, namely the cases P34F47 and P34F70, for which the highest laser power was used. For these cases, the simulated probe temperature would rise to well above \SI{100}{\celsius}, while the measured probe temperature would reach a plateau below \SI{100}{\celsius}. Therefore, in \cite{Huebner} the authors suspected that the missing vaporization model was the reason for this discrepancy. Now, we test this hypothesis by repeating the simulations with the previously introduced modified models that now include vaporization and condensation effects. 

\subsection{Experimental Setting}
\label{subsec:data}

\begin{table}[!t]
	\begin{center}
		\begin{tabular}{|l|r|r|r|r|r|r|r|r|r|}
			\hline
			Case Label & \sw{P22F47} & \sw{P22F70} & \sw{P22F92} & 
						 \sw{P28F47} & \sw{P28F70} & \sw{P28F92} & 
						 \sw{P34F47} & \sw{P34F70} & \sw{P34F92}
			\\
			\hline
			Laser Power [W] &&&&&&&&& \\
			\; -measured $\hat{q}_{app}$  & 
			22.1 & 22.1 & 22.1 & 
			28.0 & 28.0 & 28.0 & 
			33.8 & 33.8 & 33.8
			\\
			\hline
			Coolant $\dot{V}$ [ml/min] & 
			47.2 & 69.9 & 91.7 &
			47.5 & 70.3 & 91.8 &
			47.2 & 70.4 & 92.2
			\\
			\hline
			Time [s] &&&&&&&&& \\
			\; -Laser on $t_{on}$ &
			24 & 30 & 36 &
			18 & 30 & 60 &
			18 & 24 & 48  
			\\
			\; -Laser off $t_{off}$ & 
			1266 & 1236 & 684 &
			942 & 1722 & 1098 &
			1206 & 948 & 1182
			\\
			\; -End  $t_{end}$&
			1284 & 1248 & 702 &
			954 & 1734 & 1116 &
			1218 & 972 & 1206
			\\
			\hline
			Probe Position [mm] &&&&&&&&& \\
			\; -radial $d_{r}$ &
			10.1 & 11.4 & 9.2 &
			13.5 & 13.7 & 11.1 &
			11.2 & 9.9 & 9.6
			\\
			\; -axis-direction $d_{z}$ &
			12.6 & 25.7 & 20.9 &
			21.0 & 7.5 & 10.1 &
			23.8 & 26.3 & 35.3
			\\
			\hline
		\end{tabular}
	\end{center}
	\caption{Experimental setup for nine test cases (from \cite{Huebner}).}
	\label{tbl:testCases}
\end{table}

For the validation of our models, we  use the measurements from the experiments made in \cite{Huebner}. For these, livers were obtained from pigs which had been slaughtered approximately 6 hours prior to the experiment. The temperature of the samples was room temperature at the beginning of the experiments. A laser applicator from Somatex\textsuperscript{\textregistered} Medical Technologies (Teltow, Germany) was placed into the middle of the liver sample. An optical fiber from the same company with a diffuser part of 3 cm at its tip was inserted into the applicator for delivering the laser energy from a Nd:YAG laser device (MY30; Martin Medizintechnik, Tuttlingen, Germany; wavelength \SI{1064}{\nano\meter}) to the tissue. The applicator was equipped with a cooling water circulation system to protect the optical fiber and prevent the burning of tissue in close proximity to the applicator. A temperature probe was introduced into the porcine liver and placed close to the applicator in order to generate temperature measurements for validating the models of LITT.

The setup for the nine test cases is shown in Table~\ref{tbl:testCases}. The laser was applied with different powers, namely \num{22}, \num{28}, and \num{34}~\si{\watt}, and different flow rates of the applicator cooling system. However, it is assumed that the effect of the coolant flow rate is negligible (cf. \cite{Huebner}). Furthermore, the position of the temperature probe is characterized by its radial distance $d_r$ to the applicator axis as well as its distance $d_z$ from the applicator tip, where the applicator itself is contained in the half space with $z\geq 0$. We now simulate this experiment again using the vaporization models introduced previously, and compare the results with the measurement data as well as with the results obtained by the original model which does not consider vaporization. The optical, thermal, and damage parameters used for the simulation are listed in Table~\ref{table:parameters}. They are taken from \cite{Huebner} and the references therein (cf. \cite{puccini2003simulations, roggan1995optical, giering1995review, schwarzmaier1998treatment}). For the condensation region $\Omega_\text{cond}$ we have chosen the points where the temperature was between $\temperature^- = \num{60}~\si{\celsius}$ and $\temperature^+ = \num{80}~\si{\celsius}$, as proposed in \cite{Yang}.

\subsection{The Case P34F47}

\begin{figure}[!t]
	\centering
	\includegraphics[width=0.75\textwidth]{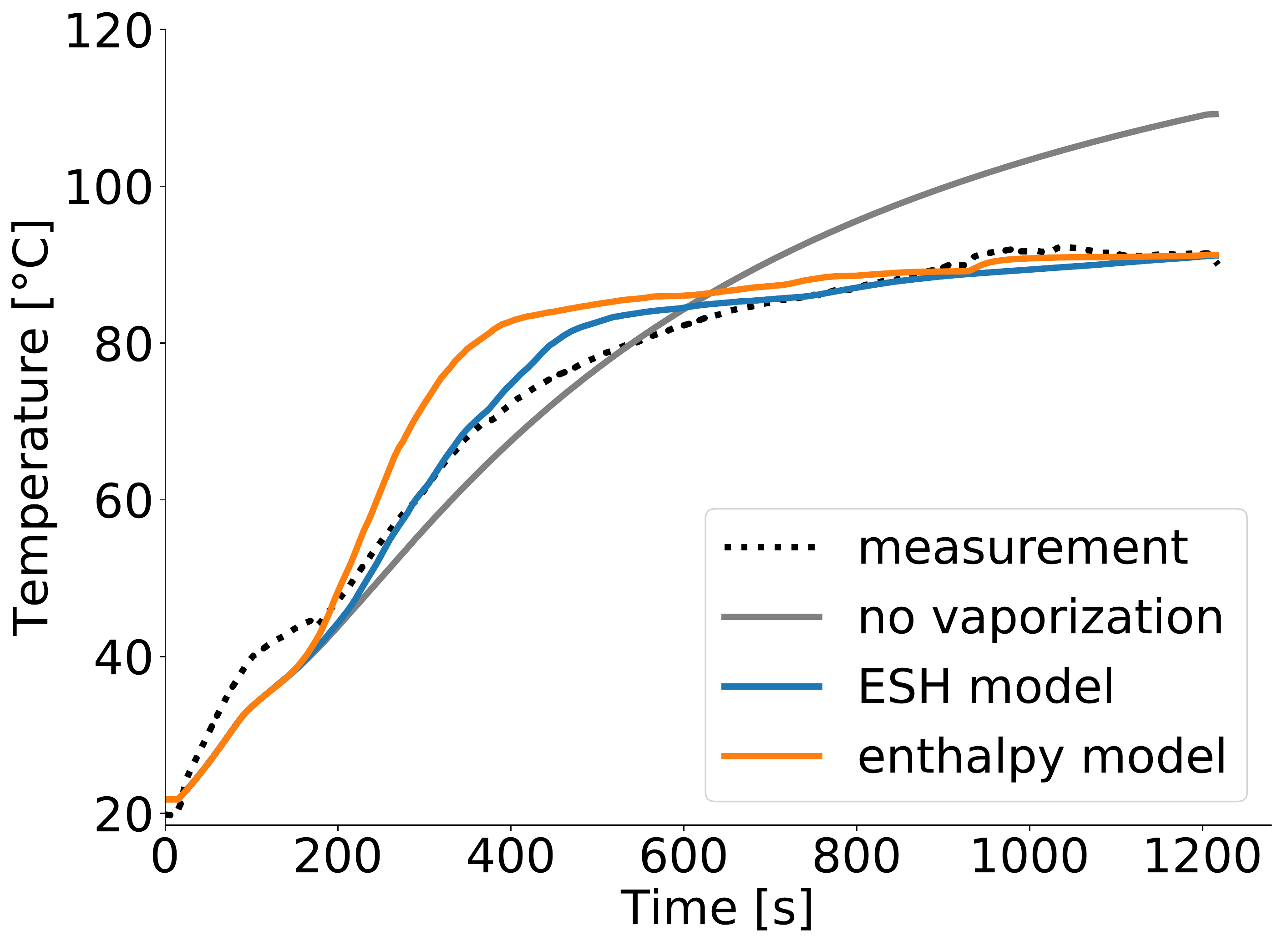}
	\caption{Comparison of the models for the case P34F47.}
	\label{fig:plot_60_80}
\end{figure}

\begin{figure}[!t]
	\centering
	\includegraphics[width=0.75\textwidth]{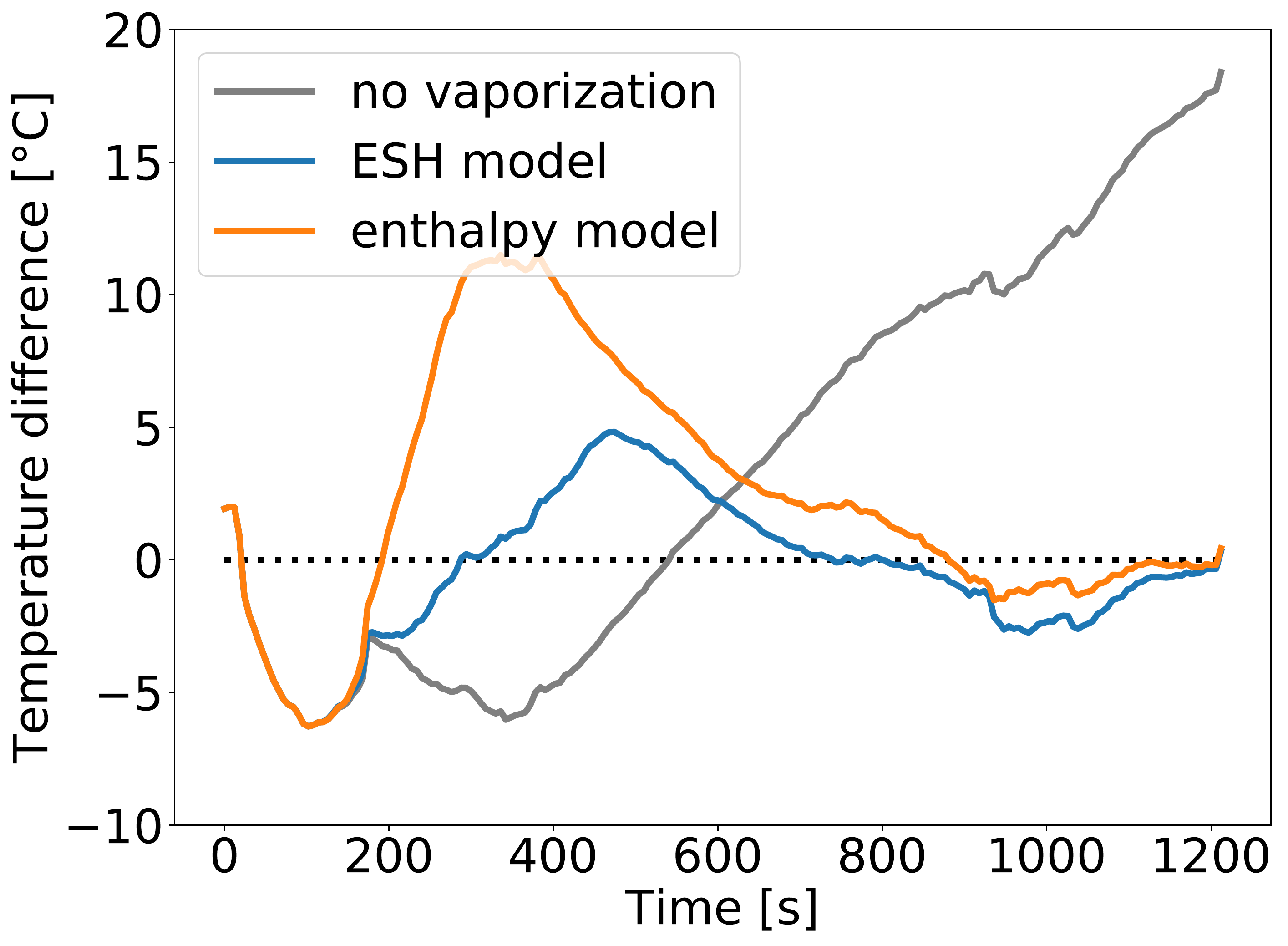}
	\caption{Difference between simulated and measured temperature for the case P34F47.}
	\label{fig:difference}
\end{figure}

Let us start the investigation of the vaporization models with the test case P34F47 of \cite{Huebner}, where a laser power of \num{34}~\si{\watt} was used. The results for this case are shown in Figure~\ref{fig:plot_60_80}, where the measurement from the temperature probe, the results for the model of \cite{Huebner} and the results for both vaporization models of Section~\ref{sec:vapo} are shown. For this specific case, the probe temperature simulated without a vaporization model rises well above \SI{100}{\celsius} while the measured temperature reaches a plateau below \SI{100}{\celsius} (see Figure~\ref{fig:plot_60_80}). In contrast, both vaporization models do not overestimate the temperature to the end of the treatment and predict the occurring plateau correctly. This is further visualized in Figure~\ref{fig:difference}, where the difference of the models from the measurement over the entire treatment is depicted. These results show that all models are reasonably close to the measured temperature until up to about \num{80}~\si{\celsius}. After that point the model without vaporization overestimates the temperature significantly. The models that include vaporization give considerably better results since their predicted temperatures match the measured ones more closely throughout the whole treatment.

\subsection{All Nine Test Cases}

\begin{figure}[!t]
	\centering
	\begin{subfigure}{0.3\textwidth}
		\includegraphics[width=\textwidth]{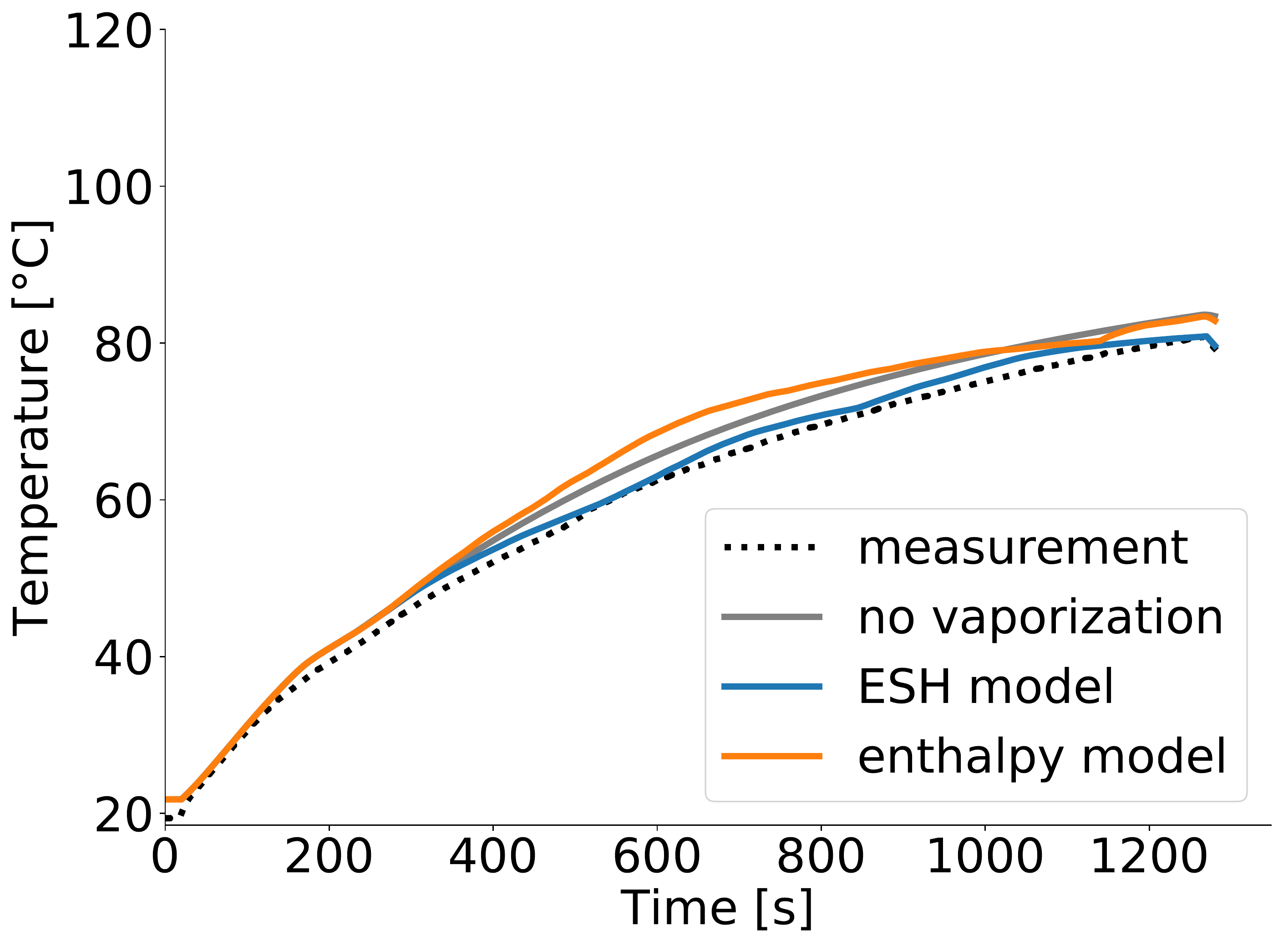}
		\caption*{P22F47}
	\end{subfigure}
	\begin{subfigure}{0.3\textwidth}
		\includegraphics[width=\textwidth]{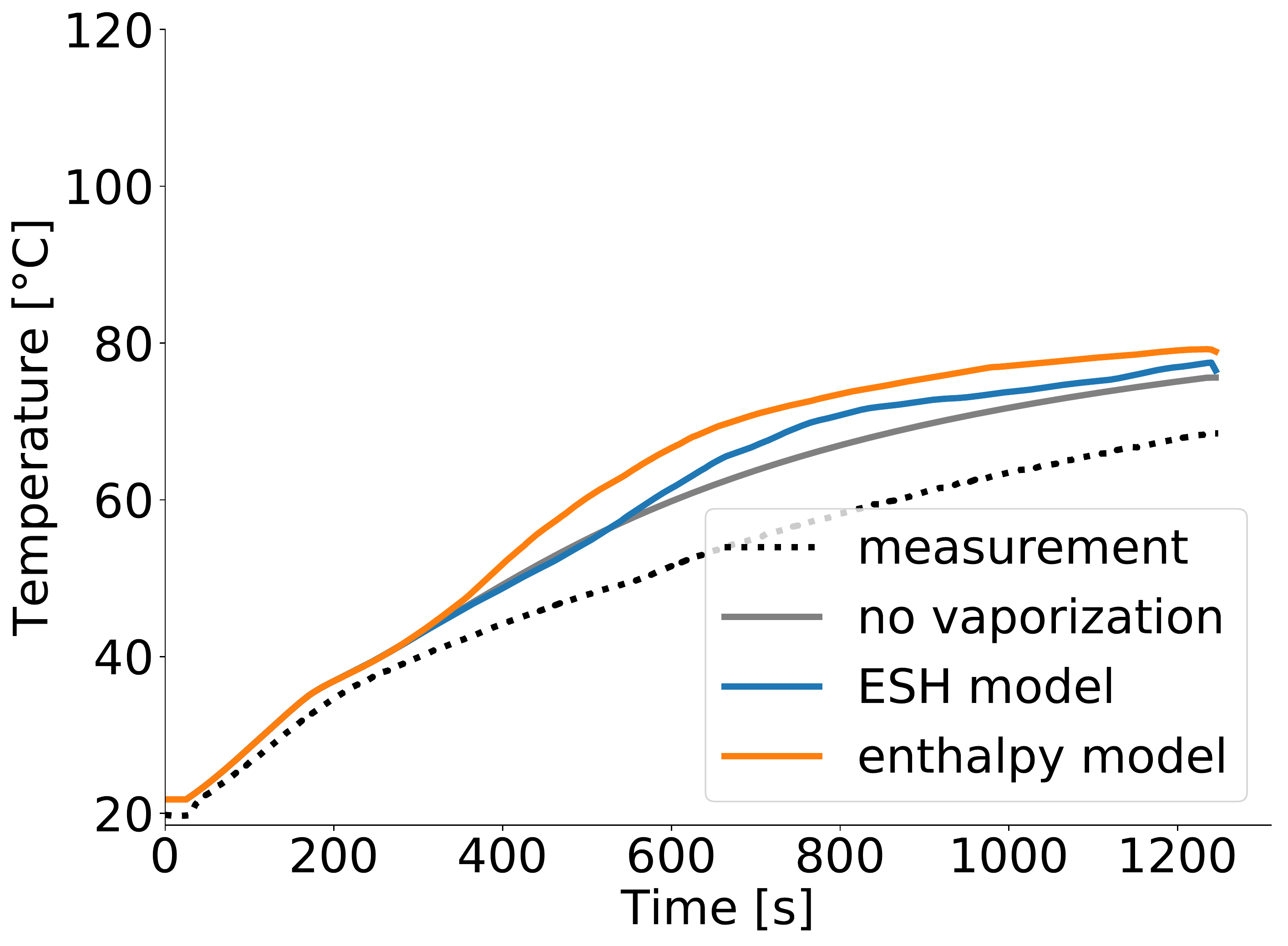}
		\caption*{P22F70}
	\end{subfigure}
	\begin{subfigure}{0.3\textwidth}
		\includegraphics[width=\textwidth]{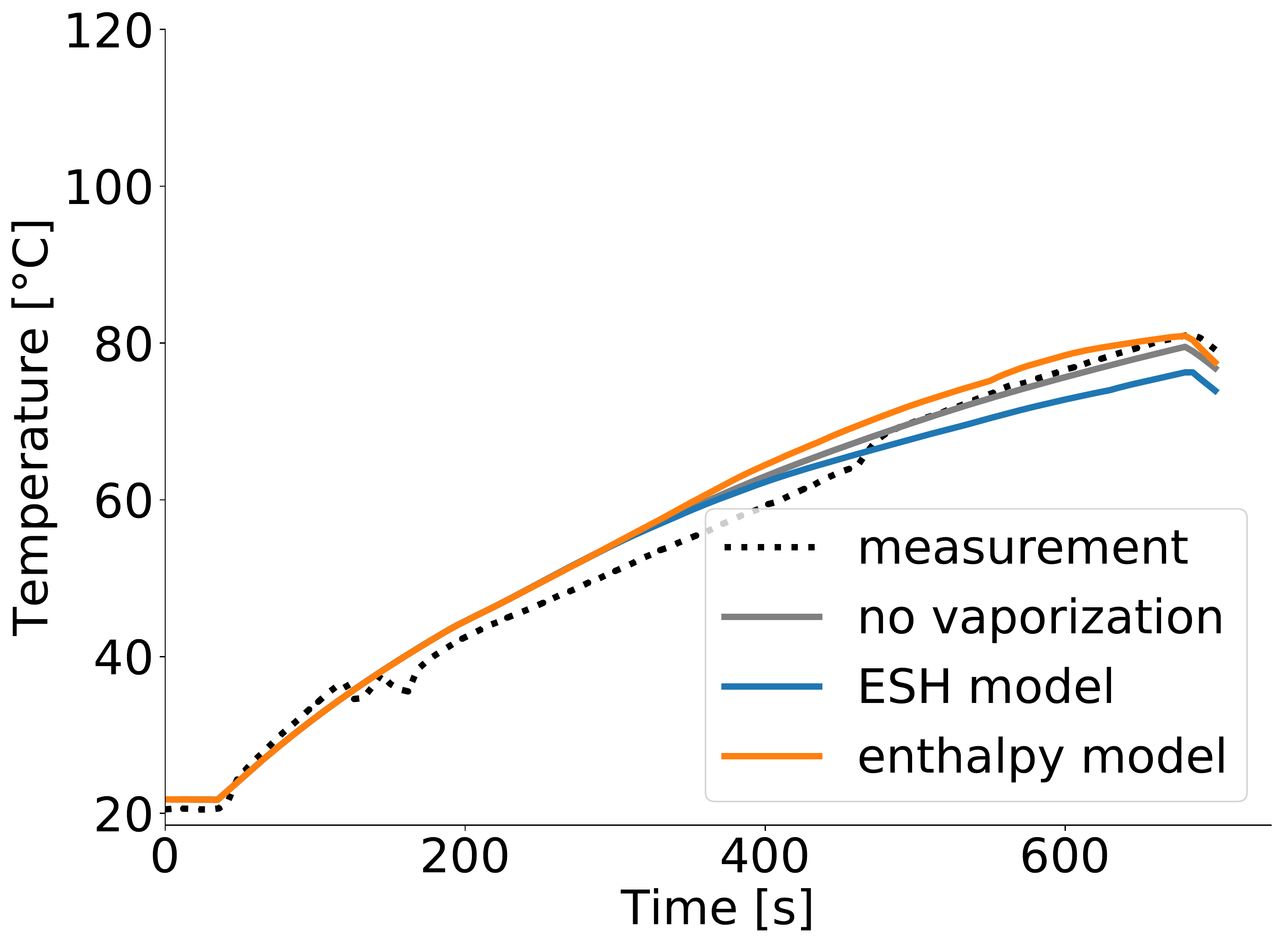}
		\caption*{P22F92}
	\end{subfigure}
	\\
	\begin{subfigure}{0.3\textwidth}
		\includegraphics[width=\textwidth]{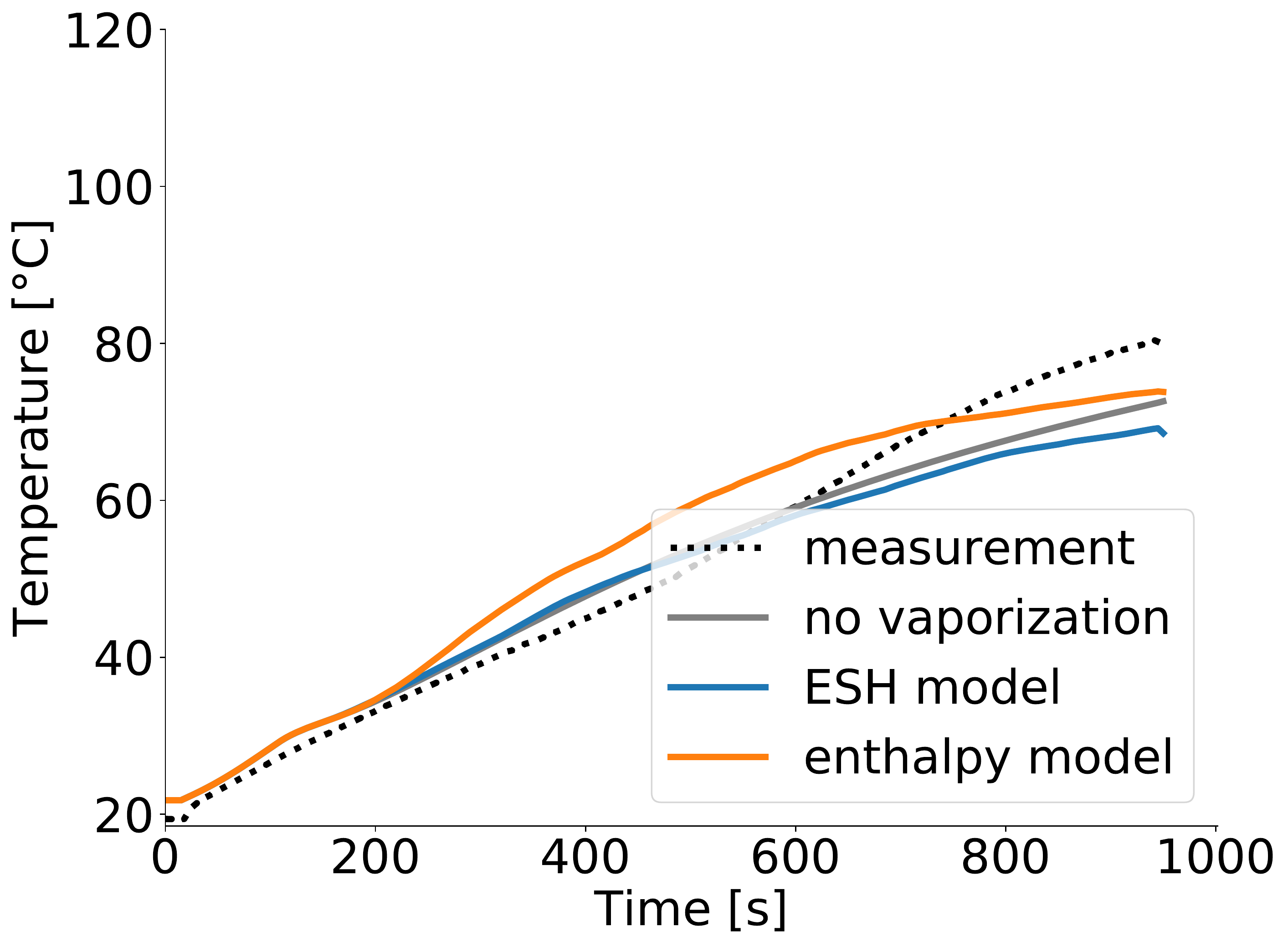}
		\caption*{P28F47}
	\end{subfigure}
	\begin{subfigure}{0.3\textwidth}
		\includegraphics[width=\textwidth]{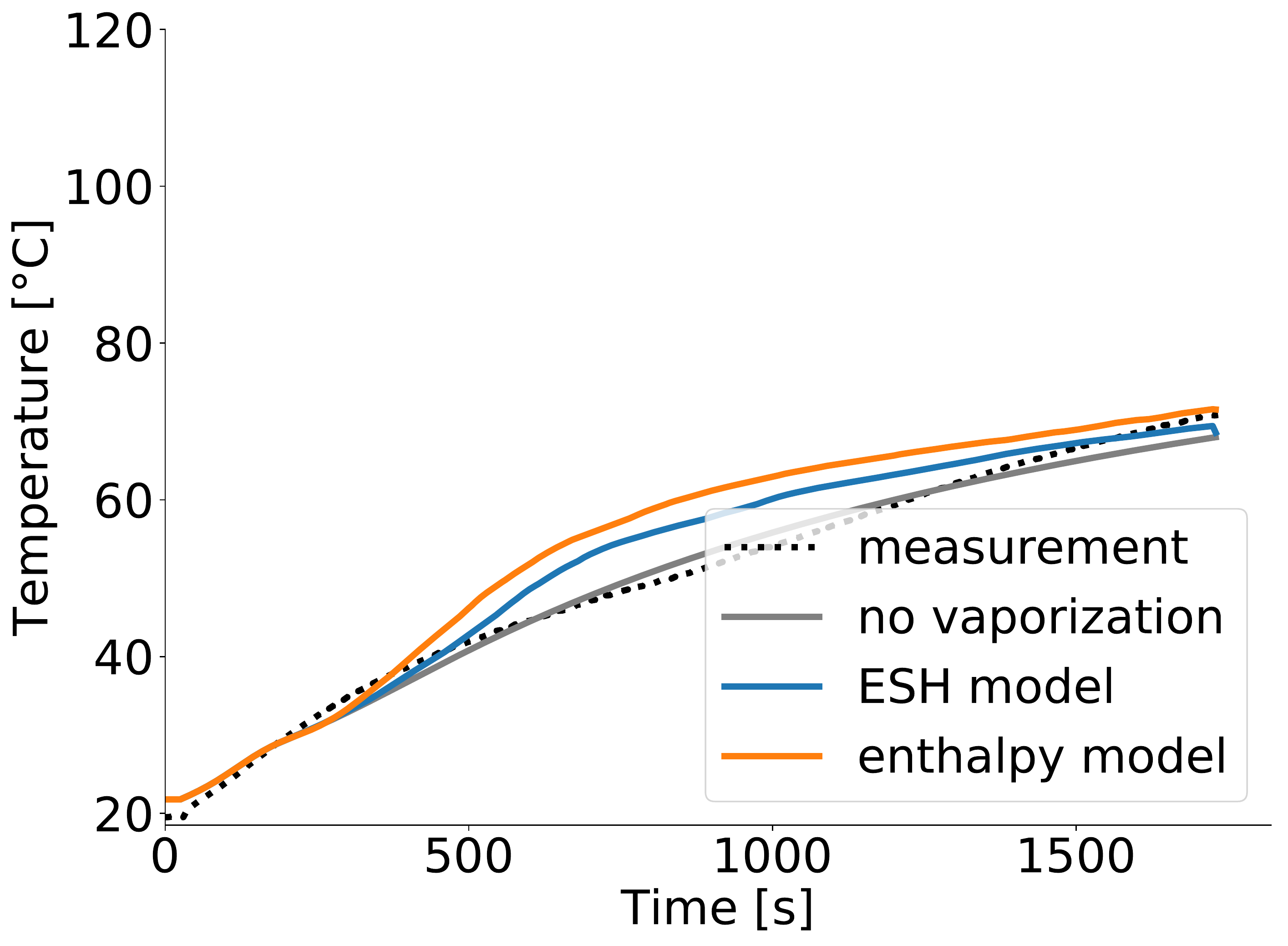}
		\caption*{P28F70}
	\end{subfigure}
	\begin{subfigure}{0.3\textwidth}
		\includegraphics[width=\textwidth]{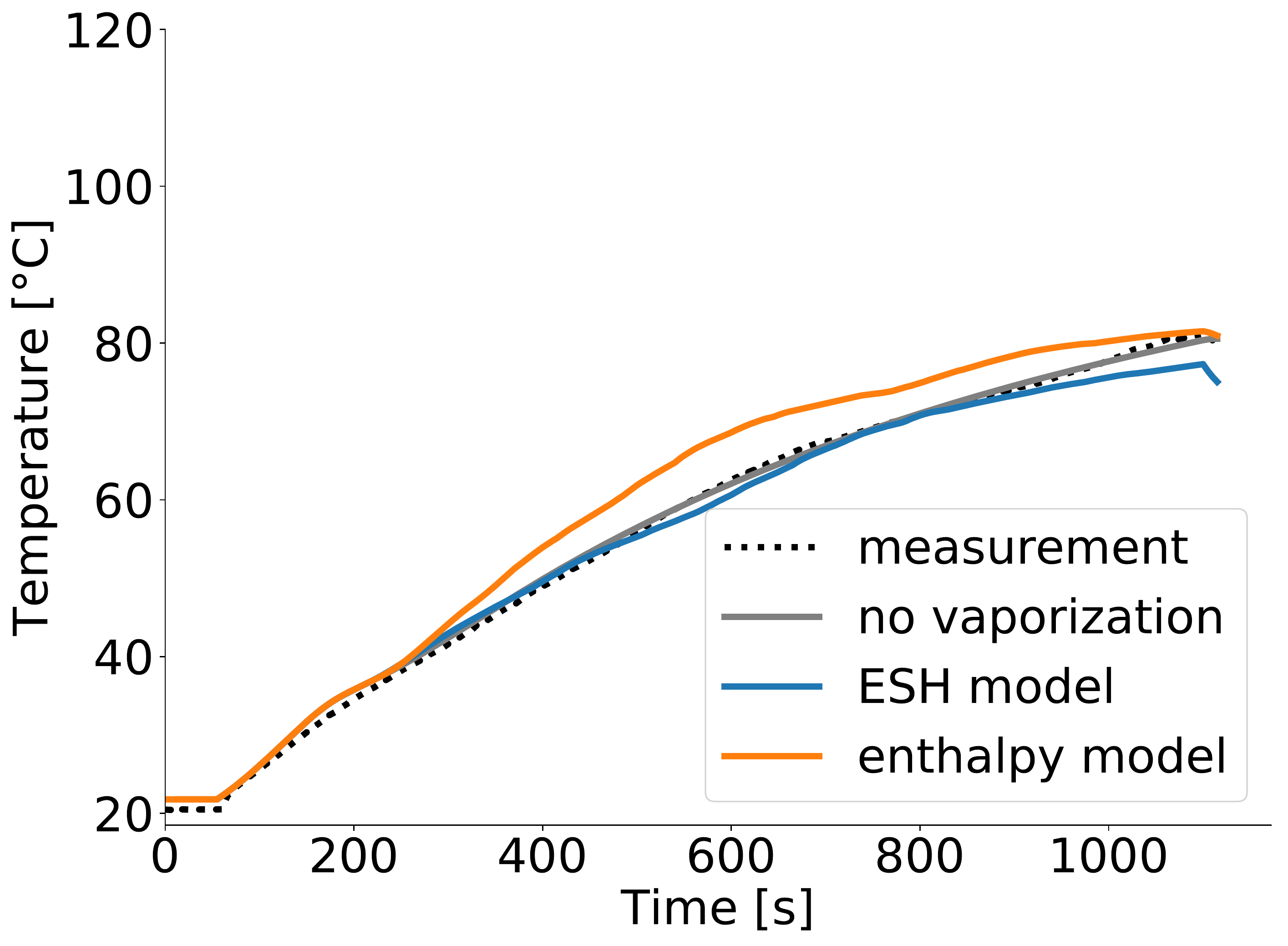}
		\caption*{P28F92}
	\end{subfigure}
	\\
	\begin{subfigure}{0.3\textwidth}
		\includegraphics[width=\textwidth]{img_60_80_K3030}
		\caption*{P34F47}
	\end{subfigure}
	\begin{subfigure}{0.3\textwidth}
		\includegraphics[width=\textwidth]{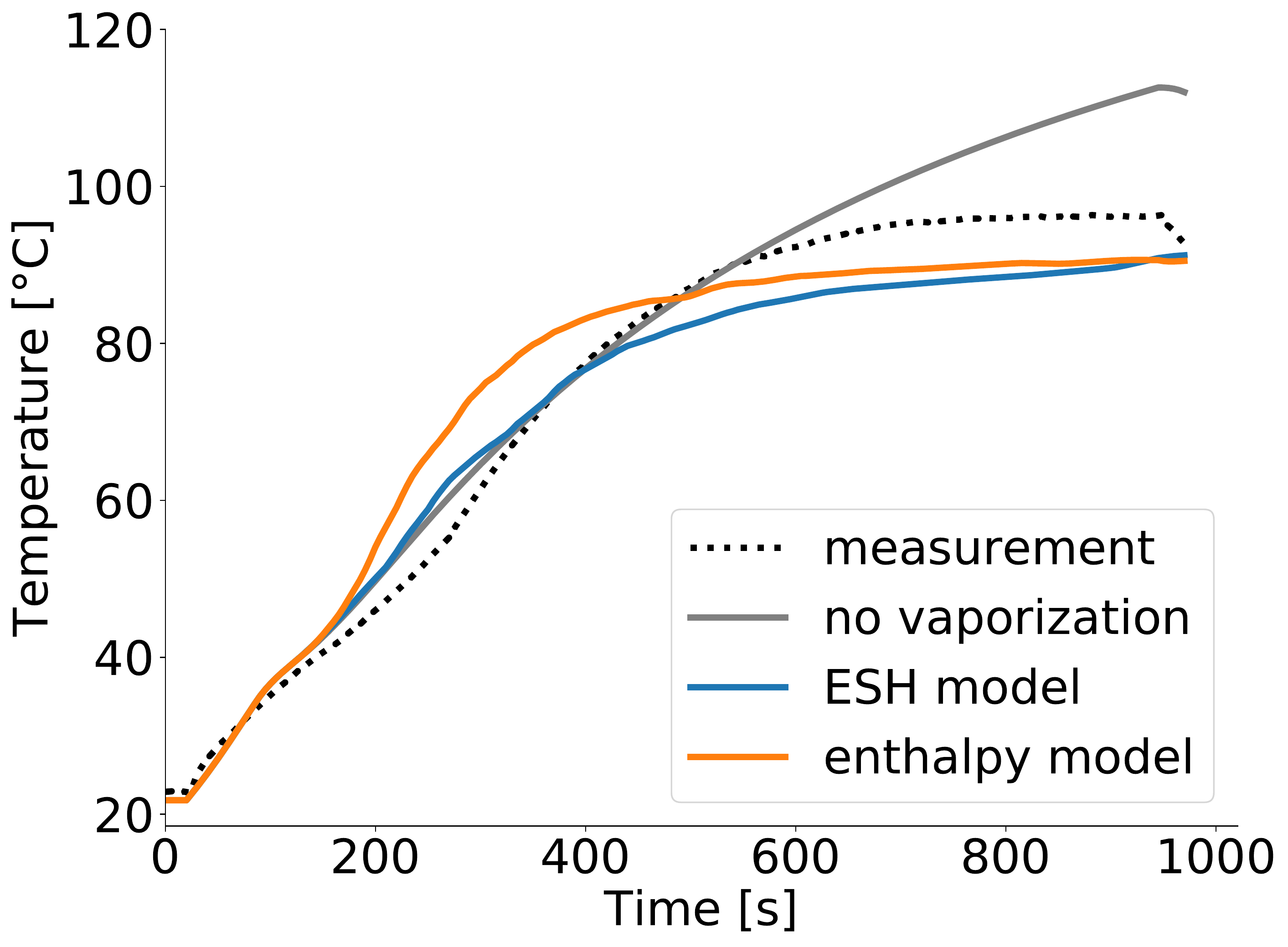}
		\caption*{P34F70}
	\end{subfigure}
	\begin{subfigure}{0.3\textwidth}
		\includegraphics[width=\textwidth]{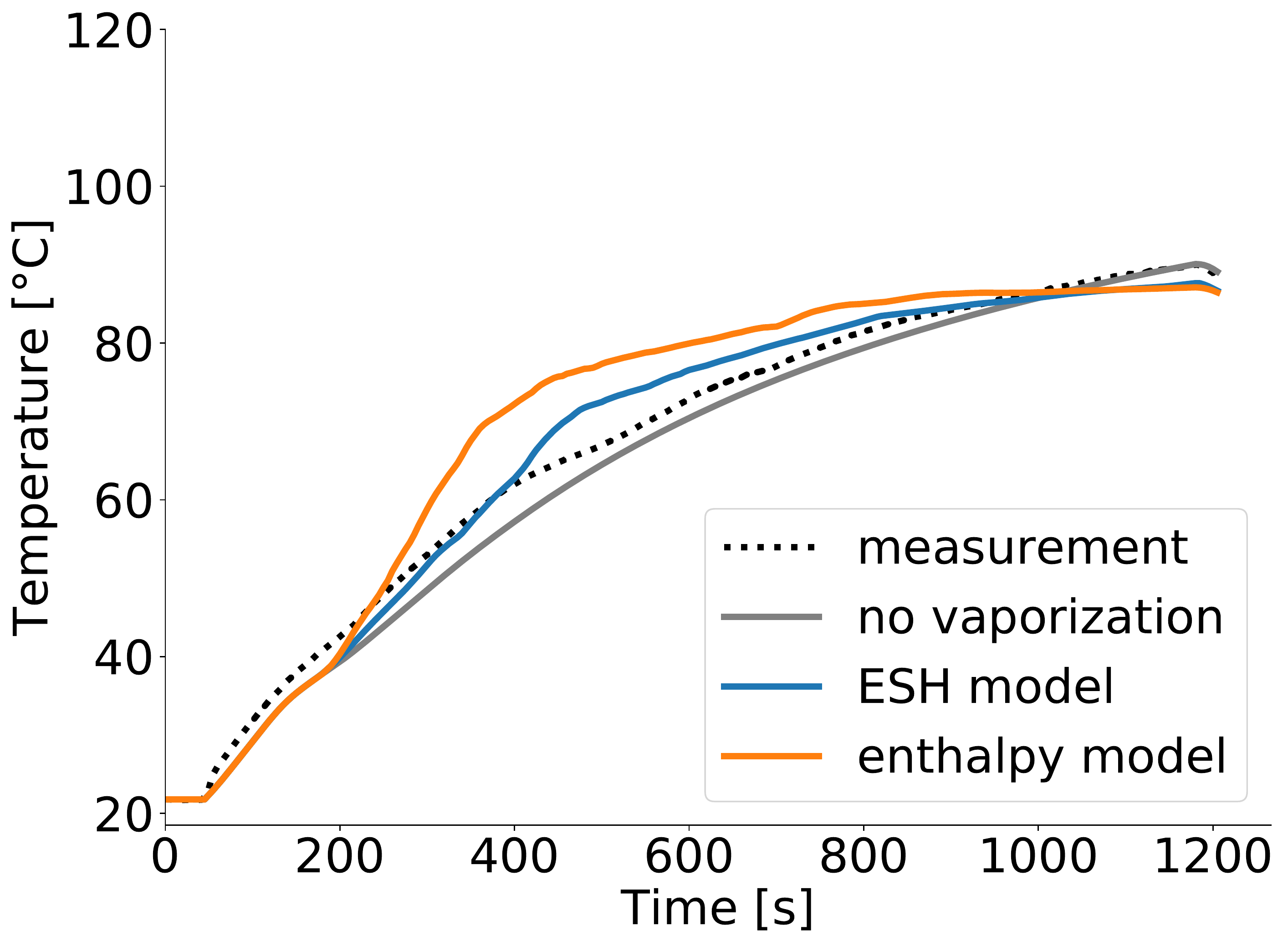}
		\caption*{P34F92}
	\end{subfigure}
	\caption{Comparison of the models with temperature measurements.}
	\label{fig:comparison9}
\end{figure}

\begin{figure}[!t]
	\centering
	\begin{subfigure}{0.3\textwidth}
		\includegraphics[width=\textwidth]{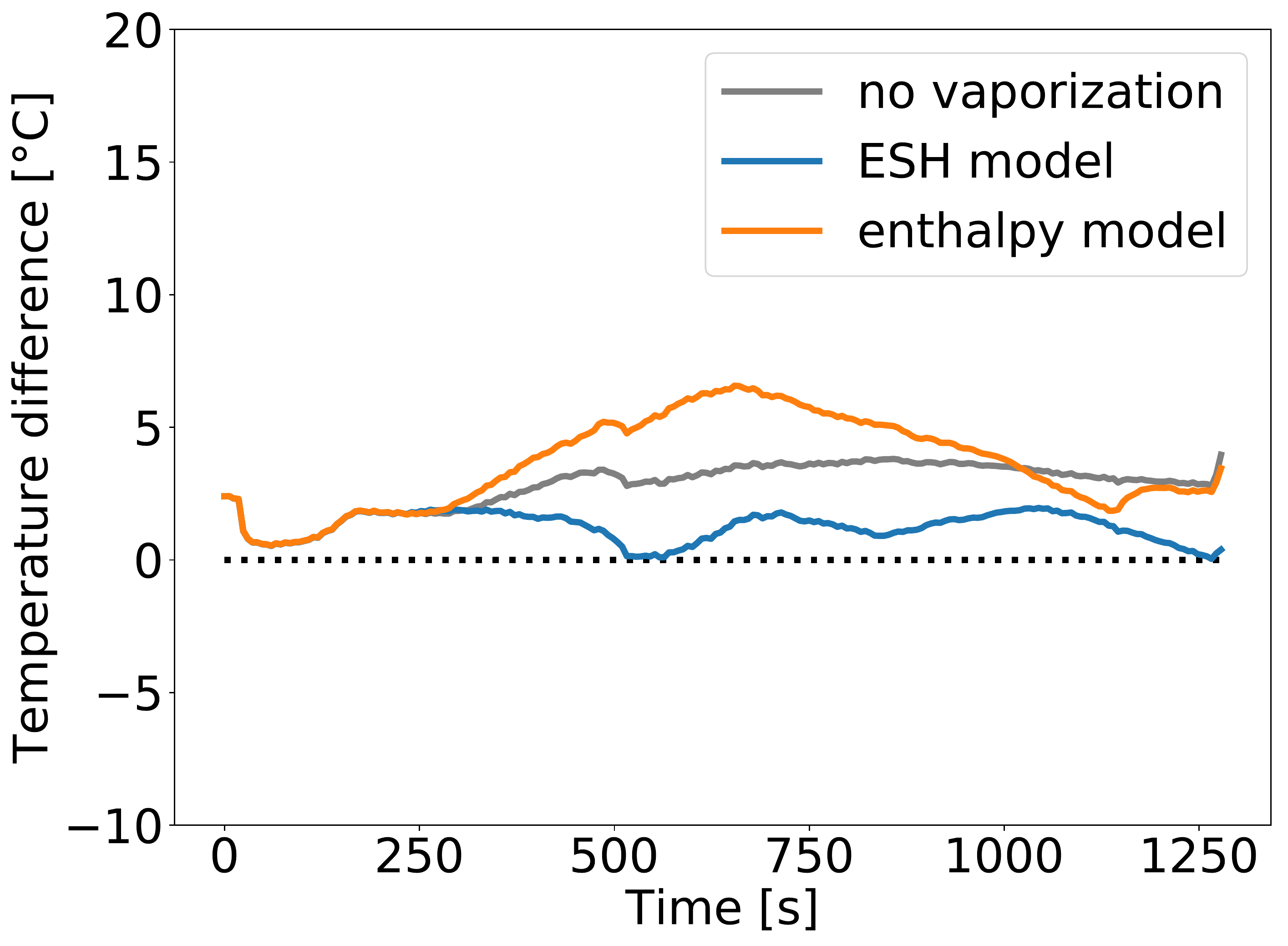}
		\caption*{P22F47}
	\end{subfigure}
	\begin{subfigure}{0.3\textwidth}
		\includegraphics[width=\textwidth]{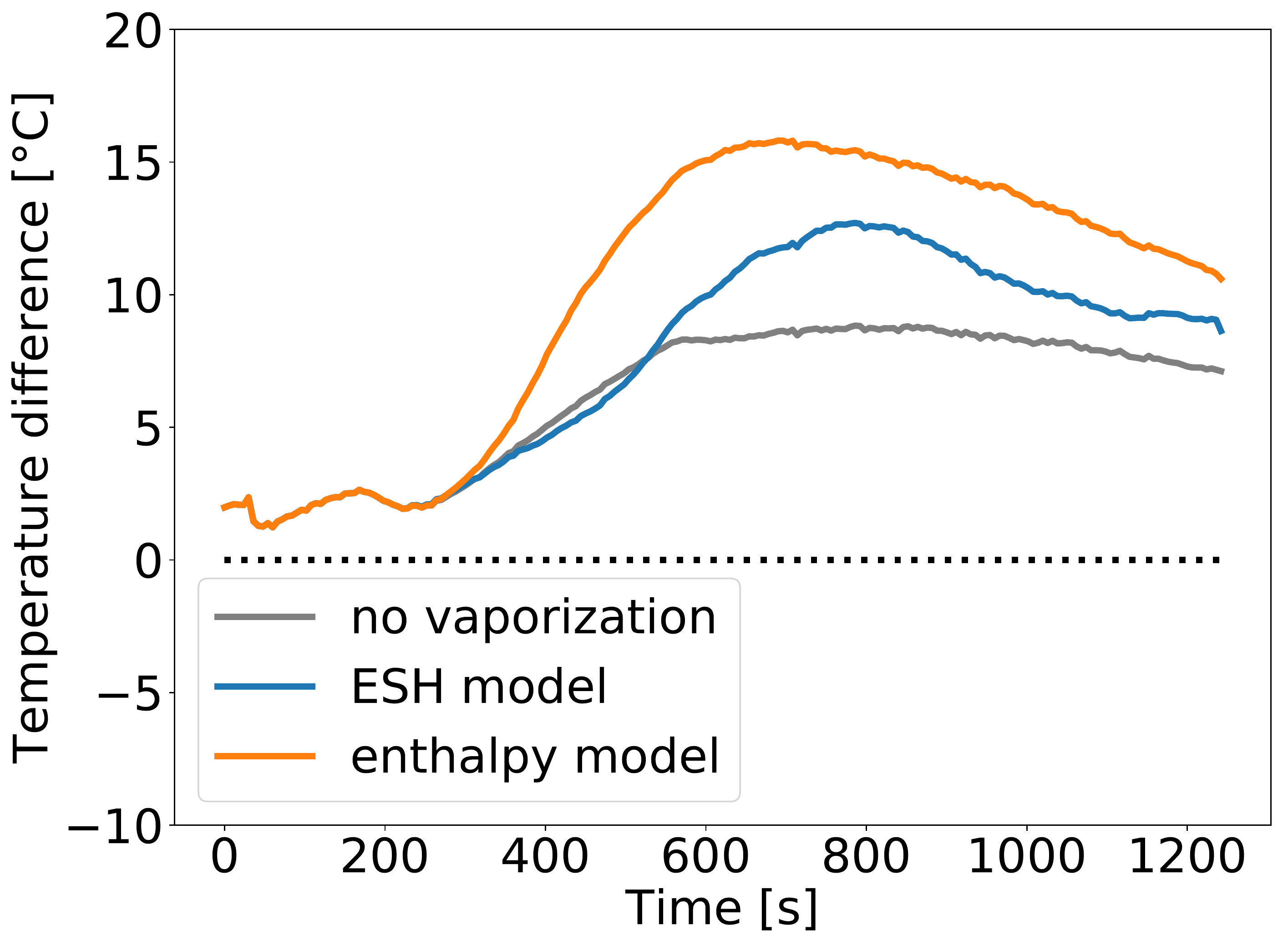}
		\caption*{P22F70}
	\end{subfigure}
	\begin{subfigure}{0.3\textwidth}
		\includegraphics[width=\textwidth]{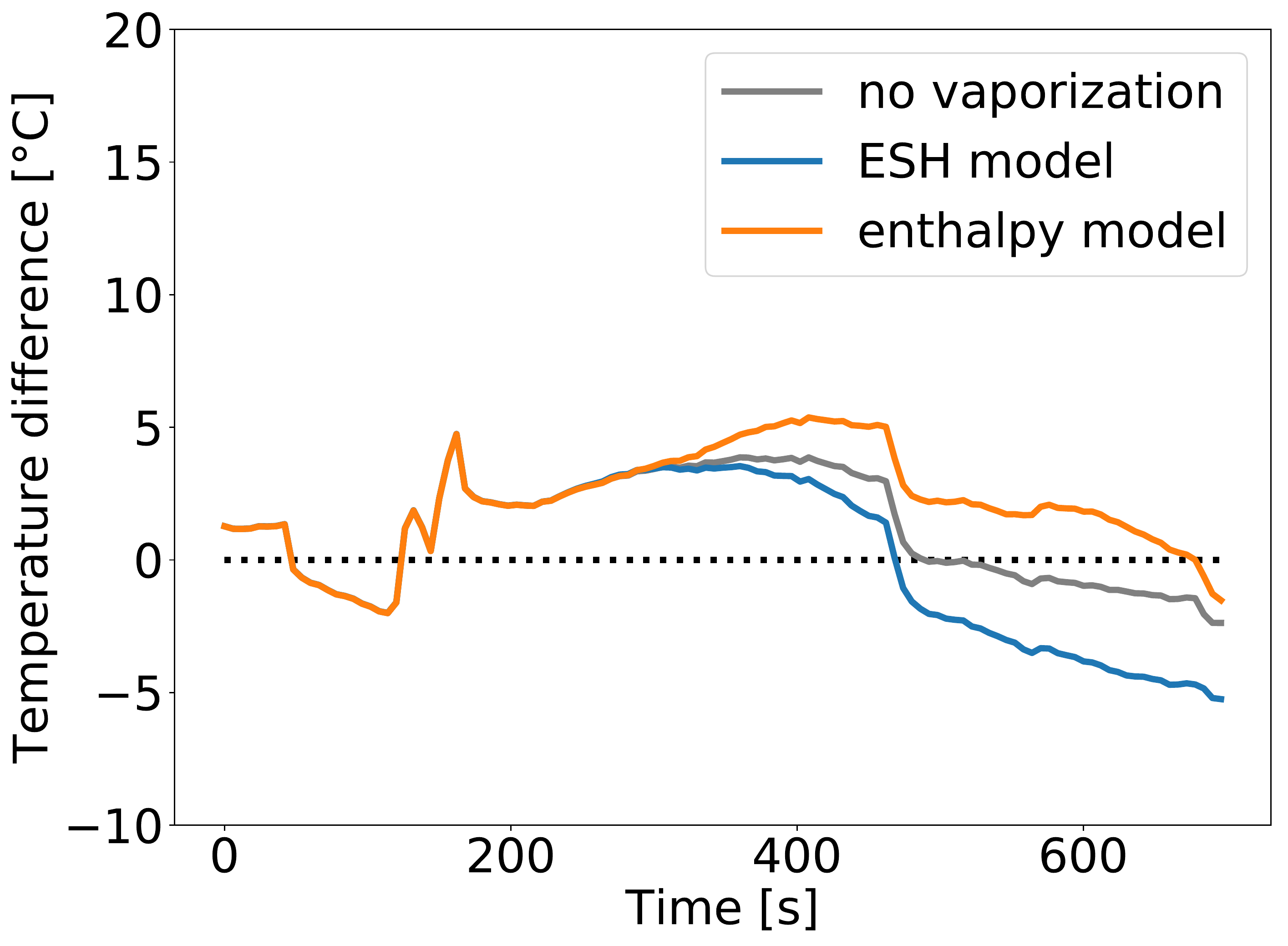}
		\caption*{P22F92}
	\end{subfigure}
	\\
	\begin{subfigure}{0.3\textwidth}
		\includegraphics[width=\textwidth]{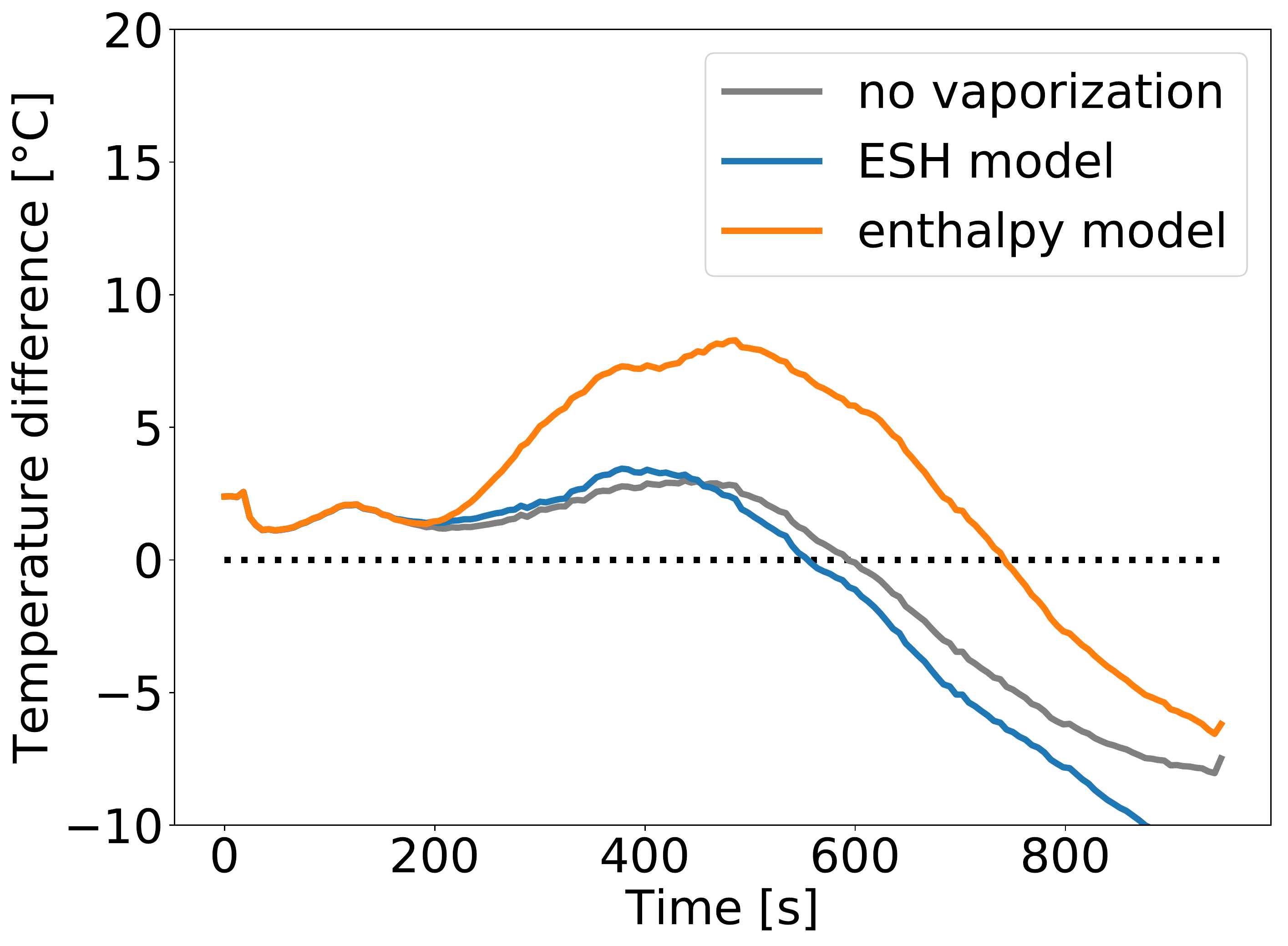}
		\caption*{P28F47}
	\end{subfigure}
	\begin{subfigure}{0.3\textwidth}
		\includegraphics[width=\textwidth]{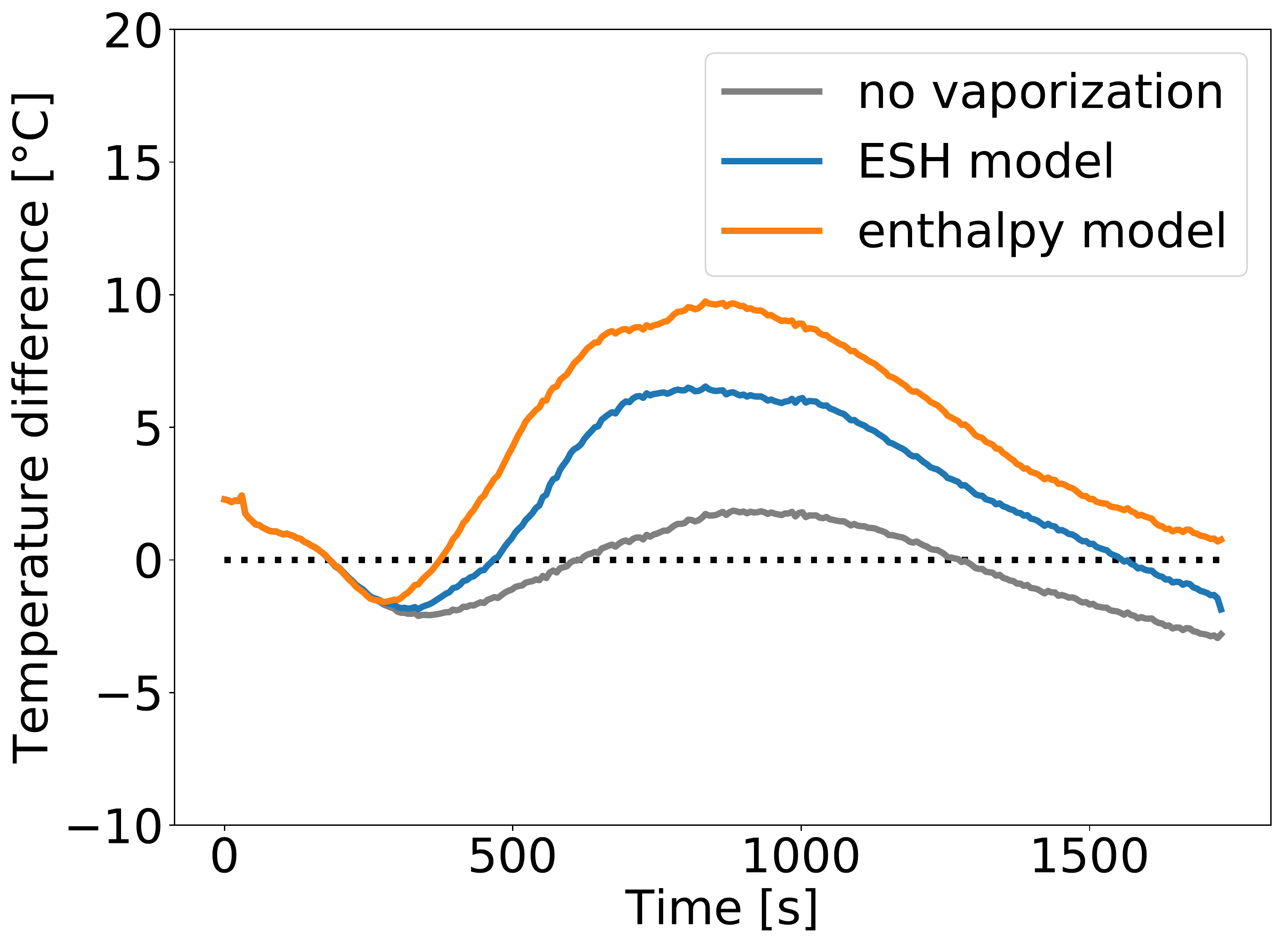}
		\caption*{P28F70}
	\end{subfigure}
	\begin{subfigure}{0.3\textwidth}
		\includegraphics[width=\textwidth]{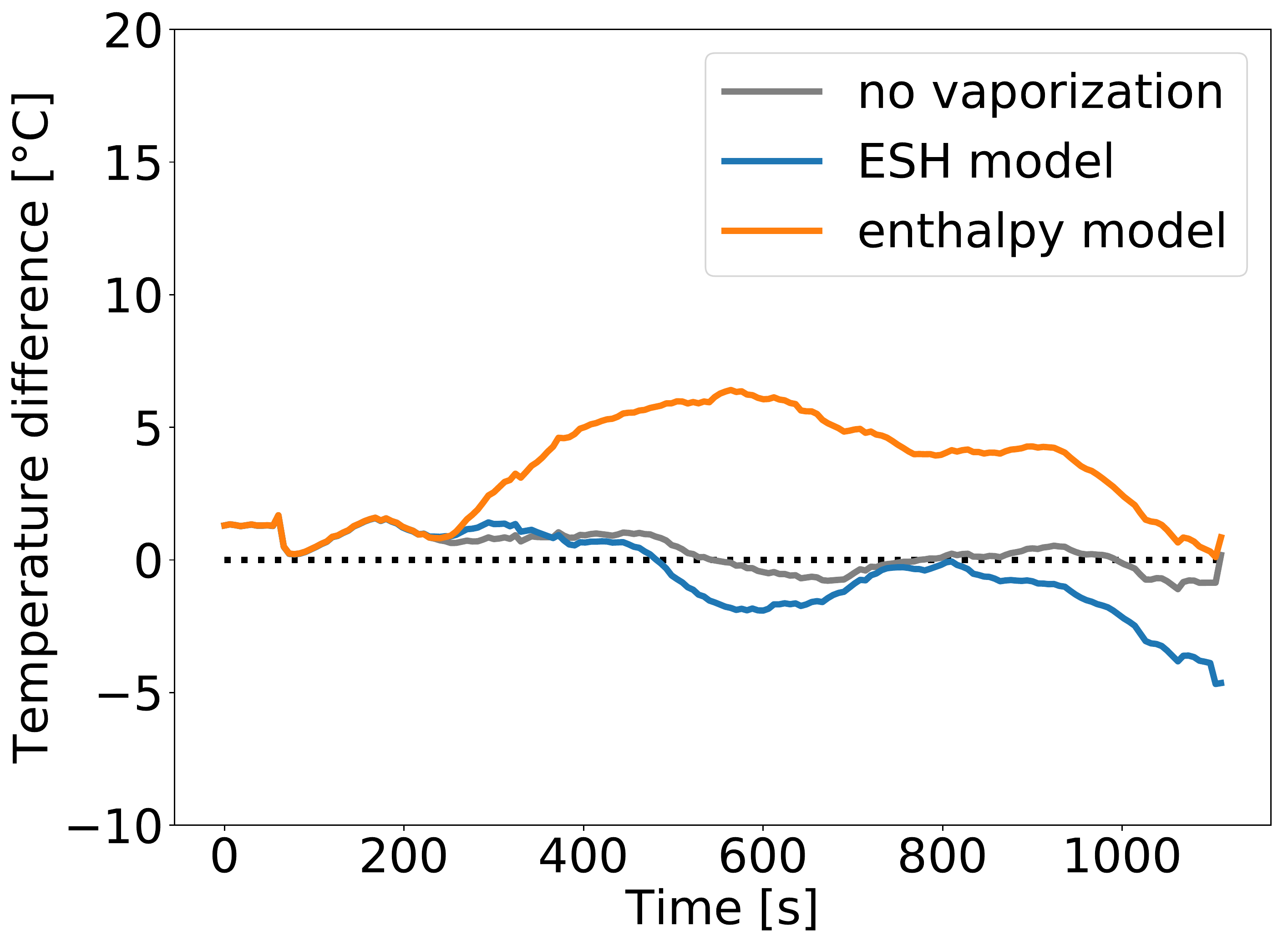}
		\caption*{P28F92}
	\end{subfigure}
	\\
	\begin{subfigure}{0.3\textwidth}
		\includegraphics[width=\textwidth]{img_60_80_diff_K3030}
		\caption*{P34F47}
	\end{subfigure}
	\begin{subfigure}{0.3\textwidth}
		\includegraphics[width=\textwidth]{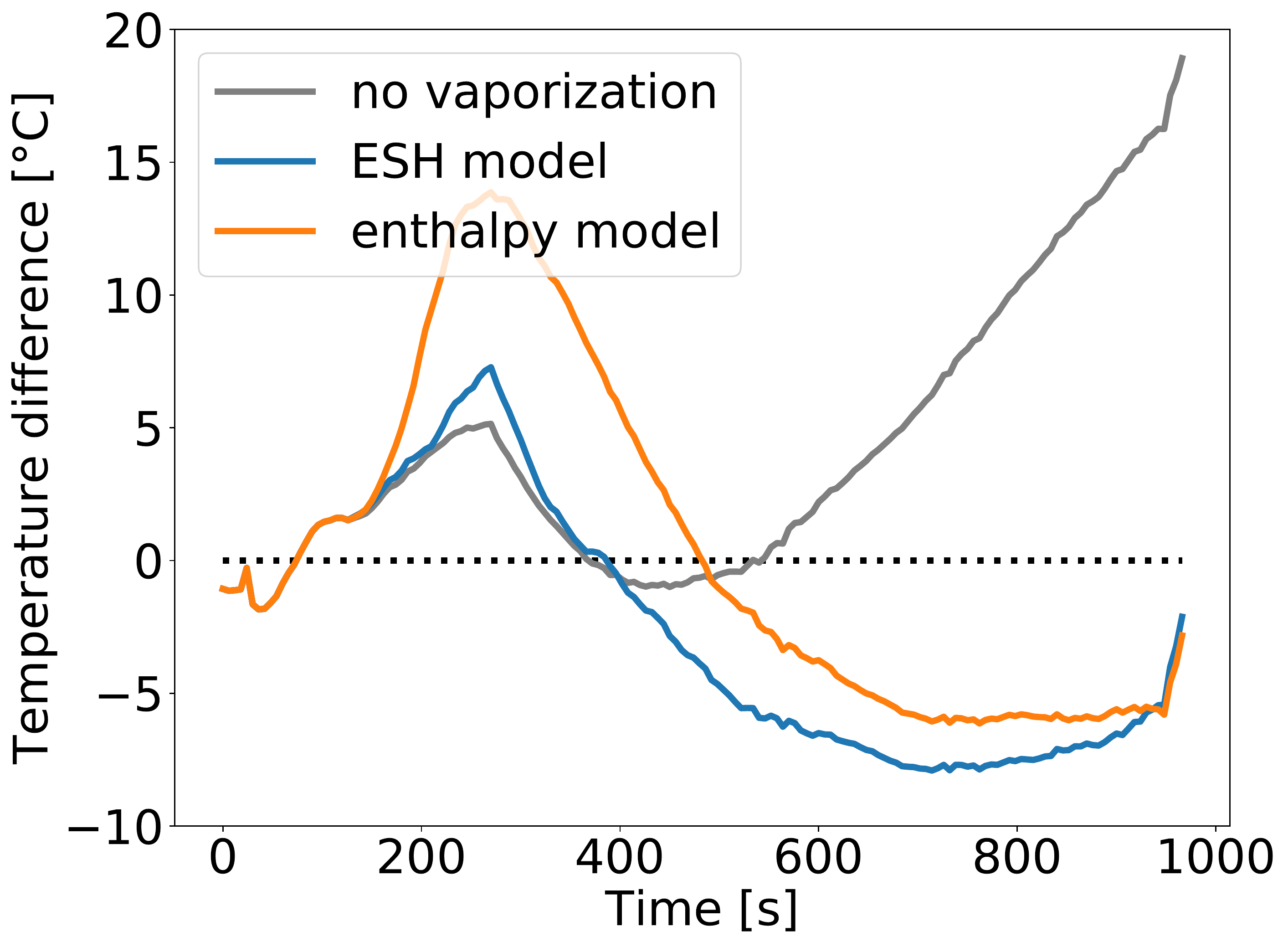}
		\caption*{P34F70}
	\end{subfigure}
	\begin{subfigure}{0.3\textwidth}
		\includegraphics[width=\textwidth]{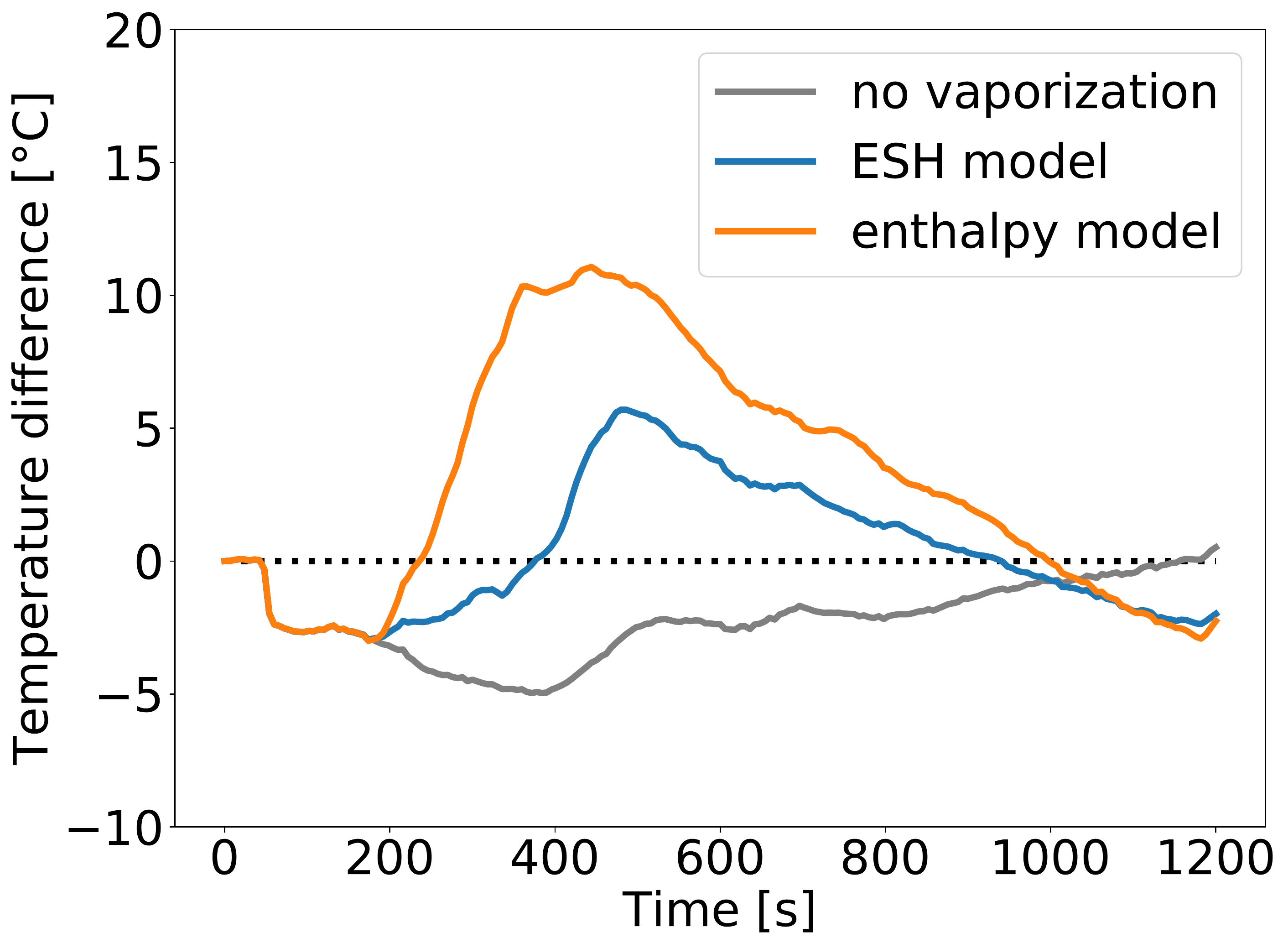}
		\caption*{P34F92}
	\end{subfigure}
	\caption{Difference between simulated and measured temperature.}
	\label{fig:difference9}
\end{figure}

After investigating the vaporization models in the context of the previous test case, where the original model without vaporization of \cite{Huebner, Fasano} failed to predict the correct temperatures, we now investigate the other test cases from the study of \cite{Huebner}. The corresponding results are shown in Figure~\ref{fig:comparison9}, where the measured and simulated temperature at the probe is shown, and Figure~\ref{fig:difference9}, which visualizes the difference of the simulated temperatures to the measurement. In general, the vaporization models are good in estimating the final temperature of the experiment. Especially for the cases P34F47 and P34F70, which could not be simulated accurately in \cite{Huebner}, the vaporization model performs significantly better and does not overestimate the temperature to the end of the treatment. However, during the middle of the experiment the vaporization models tend to overestimate the temperatures. This can be seen, e.g., for the cases P22F70 and P28F70 (cf. Figure~\ref{fig:comparison9}). We suspect that the simple condensation model is responsible for this discrepancy as we explain in Section~\ref{sec:discussion condensation model}.

Altogether, the ESH and the enthalpy model both show comparable but slightly different temperature curves. Especially the overestimation of the temperature during the middle of the experiment is usually higher for the enthalpy model. To compensate one could think about adjusting parameters, like the exact amount of water in the liver tissue. However, a first step should be to improve the simple condensation model.

\subsection{Discussion of the Simple Condensation Model}
\label{sec:discussion condensation model}
\begin{figure}[!t]
	\centering
	\begin{subfigure}{0.45\textwidth}
		\includegraphics[width=\textwidth]{img_60_80_K3030}
		\caption{$\temperature^- = \num{60}~\si{\celsius}$ and $\temperature^+ = \num{80}~\si{\celsius}$.}
		\label{fig:comp_1}
	\end{subfigure}
	\begin{subfigure}{0.45\textwidth}
		\includegraphics[width=\textwidth]{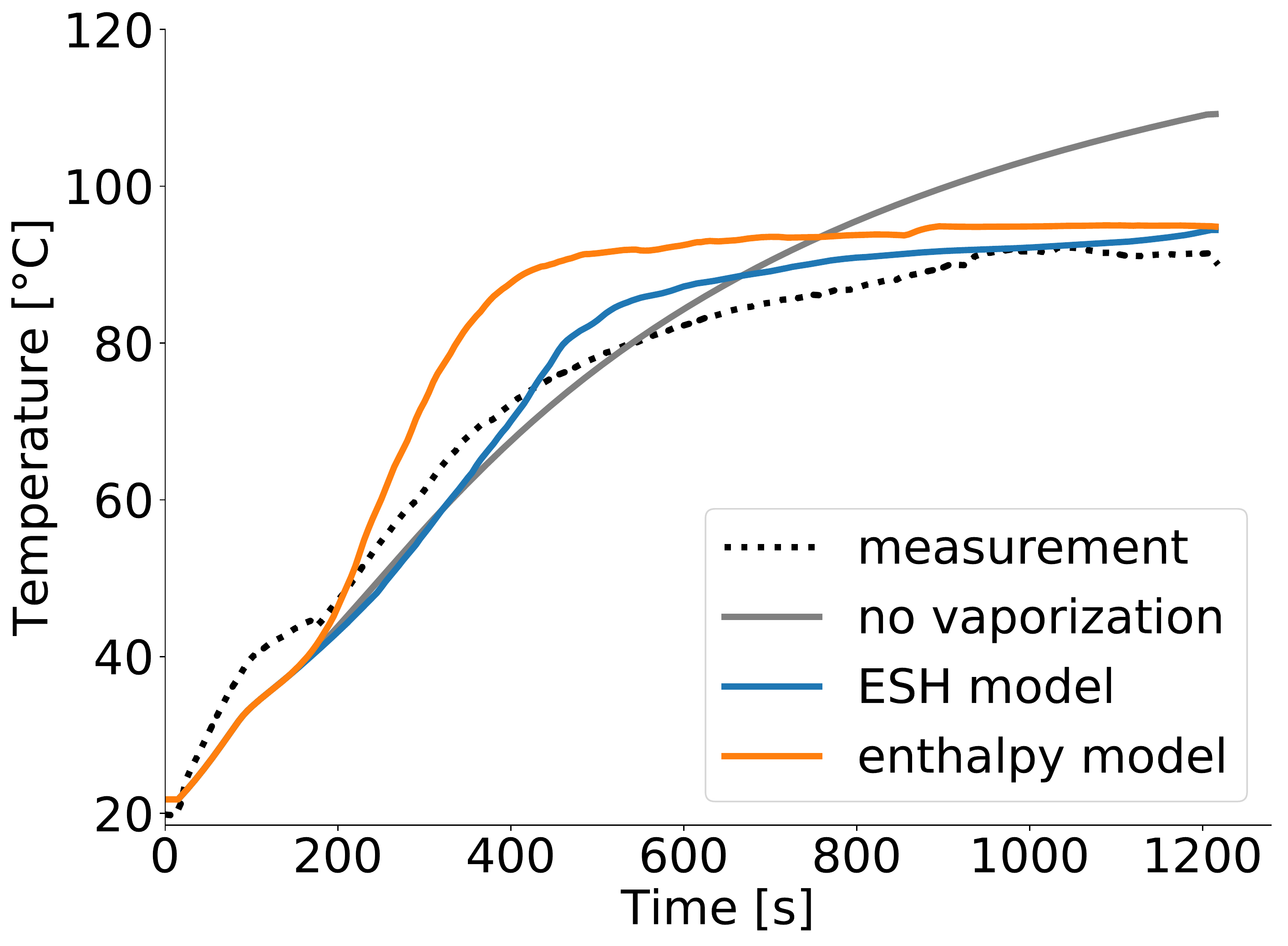}
		\caption{$\temperature^- = \num{70}~\si{\celsius}$ and $\temperature^+ = \num{90}~\si{\celsius}$.}
		\label{fig:comp_2}
	\end{subfigure}
	\caption{Sensitivities with respect to the choice of the condensation region \mbox{$\Omega_\text{cond} = \Set{x\in \Omega | \temperature^- \leq \temperature \leq \temperature^+}$}.}
	\label{fig:plot_70_90}
\end{figure}

The consideration behind the simple condensation model described in Sections~\ref{sec:condensation esh} and \ref{sec:condensation enthalpy} is solely to preserve the conservation of energy. Therefore, all the water which was vaporized at a certain time is assumed to instantaneously condensate in the condensation region $\Omega_\text{cond} = \Set{x\in \Omega | \temperature^- \leq \temperature \leq \temperature^+}$. This consideration is strictly global and does not involve any form of transport mechanism for the vapor. Hence, it is possible that vapor which was generated in one region can instantaneously condensate in another region. Through this mechanism temperature can be shifted from one region to another without any delay. This effect is possibly the reason for the overestimated temperatures during the middle of the experiment. This can be seen, e.g., for the case P28F70 (cf. Figure~\ref{fig:difference9}), where all simulated temperatures are the same until about \SI{400}{\second} into the experiment. At that point the simulated temperatures rise much faster for the models that include vaporization than for the one without it. We suspect that at this point of the experiment, vaporization occurs in tissue close to the applicator. Due to the instantaneous transport mechanism of the simplified condensation model heat is then added to regions further away from the applicator, where the applicator is placed. This results in the non-physical temperature increase that can be seen in this case. Additionally the simple condensation model is also rather sensitive with respect to the choice of the condensation region as can be seen in Figure~\ref{fig:plot_70_90}, where the temperature at the probe for the case P34F47 is shown for the condensation region given by $\temperature^- = \SI{60}{\celsius}$ and $\temperature^+ = \SI{80}{\celsius}$ in Figure~\ref{fig:comp_1} as well as for $\temperature^- = \SI{70}{\celsius}$ and $\temperature^+ = \SI{90}{\celsius}$ in Figure~\ref{fig:comp_2}.

To resolve this issue, the transport of vapor within the tissue must be taken into account. This could for example be done by adding an additional diffusion equation similar to the bio-heat equation to the state system. Therefore, an effective diffusion coefficient for the vapor must be known or estimated from measurements. Alternatively, a more complex solution would be to model the tissue as porous medium and to use a pressure based formulation for the vapor transport.

\section{Conclusion}
\label{sec:conclusion}
LITT is a minimally-invasive method in the field of interventional oncology used for treating liver cancer. Mathematical modeling and computer simulation are important features for treatment planning and simulating the necrosis of the tissue. The numerical simulation is in good agreement with temperature measurements for ex-vivo porcine liver. In particular, the incorporation of vaporization of water in liver tissue improves the simulation. However, our study suggests that the simple condensation model should be refined. Due to its global nature, this model allows for an undelayed flow of temperature from a hot region to a colder one. This is probably the reason for the overestimated temperatures during the middle of the experiments. An additional diffusion equation for the vapor with an effective or estimated diffusion coefficient could resolve this problem. In order to use simulations for the monitoring and assistance during the treatment of humans it is important to model the blood perfusion, because blood vessels have a significant cooling effect. One approach can be to identify the blood perfusion rate from MR thermometry during the beginning of the treatment and use this information to correctly simulate the remaining treatment (cf. \cite{andres2019identification}).


\begin{backmatter}
\section*{Availability of data and material}
Not applicable.

\section*{Competing interests}
The authors declare that they have no competing interests.
  
\section*{Funding}
This work was funded by the Federal Ministry of Education and Research of Germany in the framework of the project \textit{proMT: Prognostische modellbasierte online MR-Thermometrie bei minimalinvasiver Thermoablation zur Behandlung von Lebertumoren} (Förderkennzeichen: 05M16AMA). 

\section*{Author's contributions}
All authors read and approved the final manuscript.

\section*{Acknowledgments}
Not applicable.




\newcommand{\BMCxmlcomment}[1]{}

\BMCxmlcomment{

<refgrp>

<bibl id="B1">
  <title><p>Laser-induced interstitial thermotherapy</p></title>
  <aug>
    <au><snm>M{\"u}ller</snm><fnm>GJ</fnm></au>
    <au><snm>Roggan</snm><fnm>A</fnm></au>
  </aug>
  <publisher>Bellingham, WA: SPIE Press</publisher>
  <pubdate>1995</pubdate>
</bibl>

<bibl id="B2">
  <title><p>On a mathematical model for laser-induced thermotherapy</p></title>
  <aug>
    <au><snm>Fasano</snm><fnm>A</fnm></au>
    <au><snm>Hömberg</snm><fnm>D</fnm></au>
    <au><snm>Naumov</snm><fnm>D</fnm></au>
  </aug>
  <source>Applied Mathematical Modelling</source>
  <pubdate>2010</pubdate>
  <volume>34</volume>
  <issue>12</issue>
  <fpage>3831</fpage>
  <lpage>3840</lpage>
</bibl>

<bibl id="B3">
  <title><p>A finite element method model to simulate laser interstitial thermo
  therapy in anatomical inhomogeneous regions</p></title>
  <aug>
    <au><snm>Mohammed</snm><fnm>Y</fnm></au>
    <au><snm>Verhey</snm><fnm>JF</fnm></au>
  </aug>
  <source>Biomedical engineering online</source>
  <publisher>BioMed Central</publisher>
  <pubdate>2005</pubdate>
  <volume>4</volume>
  <issue>1</issue>
  <fpage>2</fpage>
</bibl>

<bibl id="B4">
  <title><p>Validation of a mathematical model for laser-induced thermotherapy
  in liver tissue</p></title>
  <aug>
    <au><snm>H{\"u}bner</snm><fnm>F.</fnm></au>
    <au><snm>Leith{\"a}user</snm><fnm>C.</fnm></au>
    <au><snm>Bazrafshan</snm><fnm>B.</fnm></au>
    <au><snm>Siedow</snm><fnm>N.</fnm></au>
    <au><snm>Vogl</snm><fnm>T. J.</fnm></au>
  </aug>
  <source>Lasers in Medical Science</source>
  <pubdate>2017</pubdate>
  <volume>32</volume>
  <issue>6</issue>
  <fpage>1399</fpage>
  <lpage>-1409</lpage>
</bibl>

<bibl id="B5">
  <title><p>Experimental Validation of a Mathematical Model for Laser-Induced
  Thermotherapy</p></title>
  <aug>
    <au><snm>Leith\"auser</snm><fnm>C.</fnm></au>
    <au><snm>H\"ubner</snm><fnm>F.</fnm></au>
    <au><snm>Bazrafshan</snm><fnm>B.</fnm></au>
    <au><snm>Siedow</snm><fnm>N.</fnm></au>
    <au><snm>Vogl</snm><fnm>T.J.</fnm></au>
  </aug>
  <source>European Consortium for Mathematics in Industry</source>
  <publisher>Heidelberg: Springer</publisher>
  <pubdate>2018</pubdate>
</bibl>

<bibl id="B6">
  <title><p>Expanding the Bioheat Equation to Include Tissue Internal Water
  Evaporation During Heating</p></title>
  <aug>
    <au><snm>{Yang}</snm><fnm>D.</fnm></au>
    <au><snm>{Converse}</snm><fnm>M. C.</fnm></au>
    <au><snm>{Mahvi}</snm><fnm>D. M.</fnm></au>
    <au><snm>{Webster}</snm><fnm>J. G.</fnm></au>
  </aug>
  <source>IEEE Transactions on Biomedical Engineering</source>
  <pubdate>2007</pubdate>
  <volume>54</volume>
  <issue>8</issue>
  <fpage>1382</fpage>
  <lpage>1388</lpage>
</bibl>

<bibl id="B7">
  <title><p>Analysis of Tissue and Arterial Blood Temperatures in the Resting
  Human Forearm</p></title>
  <aug>
    <au><snm>Pennes</snm><fnm>HH</fnm></au>
  </aug>
  <source>Journal of Applied Physiology</source>
  <pubdate>1948</pubdate>
  <volume>1</volume>
  <issue>2</issue>
  <fpage>93</fpage>
  <lpage>122</lpage>
</bibl>

<bibl id="B8">
  <title><p>Laser-tissue interactions</p></title>
  <aug>
    <au><snm>Niemz</snm><fnm>MH</fnm></au>
    <au><cnm>others</cnm></au>
  </aug>
  <publisher>Berlin Heidelberg: Springer</publisher>
  <pubdate>2007</pubdate>
</bibl>

<bibl id="B9">
  <title><p>Radiative Heat Transfer</p></title>
  <aug>
    <au><snm>Modest</snm><fnm>M.F.</fnm></au>
  </aug>
  <publisher>San Diego: Academic Press</publisher>
  <pubdate>2003</pubdate>
</bibl>

<bibl id="B10">
  <title><p>Experimentalphysik 1</p></title>
  <aug>
    <au><snm>Demtr\"{o}der</snm><fnm>W</fnm></au>
  </aug>
  <publisher>Berlin Heidelberg: Springer</publisher>
  <pubdate>2018</pubdate>
</bibl>

<bibl id="B11">
  <title><p>Measurement and Analysis of Tissue Temperature During Microwave
  Liver Ablation</p></title>
  <aug>
    <au><snm>{Yang}</snm><fnm>D.</fnm></au>
    <au><snm>{Converse}</snm><fnm>M. C.</fnm></au>
    <au><snm>{Mahvi}</snm><fnm>D. M.</fnm></au>
    <au><snm>{Webster}</snm><fnm>J. G.</fnm></au>
  </aug>
  <source>IEEE Transactions on Biomedical Engineering</source>
  <pubdate>2007</pubdate>
  <volume>54</volume>
  <issue>1</issue>
  <fpage>150</fpage>
  <lpage>155</lpage>
</bibl>

<bibl id="B12">
  <title><p>Gmsh: {A} 3-{D} finite element mesh generator with built-in pre-
  and post-processing facilities</p></title>
  <aug>
    <au><snm>Geuzaine</snm><fnm>C</fnm></au>
    <au><snm>Remacle</snm><fnm>JF</fnm></au>
  </aug>
  <source>Internat. J. Numer. Methods Engrg.</source>
  <pubdate>2009</pubdate>
  <volume>79</volume>
  <issue>11</issue>
  <fpage>1309</fpage>
  <lpage>-1331</lpage>
</bibl>

<bibl id="B13">
  <title><p>The FEniCS Project Version 1.5</p></title>
  <aug>
    <au><snm>Aln{\ae}s</snm><fnm>MS</fnm></au>
    <au><snm>Blechta</snm><fnm>J</fnm></au>
    <au><snm>Hake</snm><fnm>J</fnm></au>
    <au><snm>Johansson</snm><fnm>A</fnm></au>
    <au><snm>Kehlet</snm><fnm>B</fnm></au>
    <au><snm>Logg</snm><fnm>A</fnm></au>
    <au><snm>Richardson</snm><fnm>C</fnm></au>
    <au><snm>Ring</snm><fnm>J</fnm></au>
    <au><snm>Rognes</snm><fnm>ME</fnm></au>
    <au><snm>Wells</snm><fnm>GN</fnm></au>
  </aug>
  <source>Archive of Numerical Software</source>
  <pubdate>2015</pubdate>
  <volume>3</volume>
  <issue>100</issue>
</bibl>

<bibl id="B14">
  <title><p>Automated Solution of Differential Equations by the Finite Element
  Method</p></title>
  <aug>
    <au><snm>Logg</snm><fnm>A</fnm></au>
    <au><snm>Mardal</snm><fnm>KA</fnm></au>
    <au><snm>Wells</snm><fnm>GN</fnm></au>
    <au><cnm>others</cnm></au>
  </aug>
  <publisher>Heidelberg: Springer</publisher>
  <pubdate>2012</pubdate>
</bibl>

<bibl id="B15">
  <title><p>{PETS}c Users Manual</p></title>
  <aug>
    <au><snm>Balay</snm><fnm>S</fnm></au>
    <au><snm>Abhyankar</snm><fnm>S</fnm></au>
    <au><snm>Adams</snm><fnm>MF</fnm></au>
    <au><snm>Brown</snm><fnm>J</fnm></au>
    <au><snm>Brune</snm><fnm>P</fnm></au>
    <au><snm>Buschelman</snm><fnm>K</fnm></au>
    <au><snm>Dalcin</snm><fnm>L</fnm></au>
    <au><snm>Dener</snm><fnm>A</fnm></au>
    <au><snm>Eijkhout</snm><fnm>V</fnm></au>
    <au><snm>Gropp</snm><fnm>WD</fnm></au>
    <au><snm>Karpeyev</snm><fnm>D</fnm></au>
    <au><snm>Kaushik</snm><fnm>D</fnm></au>
    <au><snm>Knepley</snm><fnm>MG</fnm></au>
    <au><snm>May</snm><fnm>DA</fnm></au>
    <au><snm>McInnes</snm><fnm>LC</fnm></au>
    <au><snm>Mills</snm><fnm>RT</fnm></au>
    <au><snm>Munson</snm><fnm>T</fnm></au>
    <au><snm>Rupp</snm><fnm>K</fnm></au>
    <au><snm>Sanan</snm><fnm>P</fnm></au>
    <au><snm>Smith</snm><fnm>BF</fnm></au>
    <au><snm>Zampini</snm><fnm>S</fnm></au>
    <au><snm>Zhang</snm><fnm>H</fnm></au>
    <au><snm>Zhang</snm><fnm>H</fnm></au>
  </aug>
  <pubdate>2019</pubdate>
  <issue>ANL-95/11 - Revision 3.11</issue>
  <url>https://www.mcs.anl.gov/petsc</url>
</bibl>

<bibl id="B16">
  <title><p>Simulations of thermal tissue coagulation and their value for the
  planning and monitoring of laser-induced interstitial thermotherapy
  (LITT)</p></title>
  <aug>
    <au><snm>Puccini</snm><fnm>S</fnm></au>
    <au><snm>B\"ar</snm><fnm>NK</fnm></au>
    <au><snm>Bublat</snm><fnm>M</fnm></au>
    <au><snm>Kahn</snm><fnm>T</fnm></au>
    <au><snm>Busse</snm><fnm>H</fnm></au>
  </aug>
  <source>Magnetic Resonance in Medicine</source>
  <pubdate>2003</pubdate>
  <volume>49</volume>
  <issue>2</issue>
  <fpage>351</fpage>
  <lpage>362</lpage>
</bibl>

<bibl id="B17">
  <title><p>The optical properties of biological tissue in the near infrared
  wavelength range</p></title>
  <aug>
    <au><snm>Roggan</snm><fnm>A</fnm></au>
    <au><snm>Dorschel</snm><fnm>K</fnm></au>
    <au><snm>Minet</snm><fnm>O</fnm></au>
    <au><snm>Wolff</snm><fnm>D</fnm></au>
    <au><snm>Muller</snm><fnm>G</fnm></au>
  </aug>
  <source>Laser-induced interstitial therapy. SPIE Press, Bellingham,
  WA</source>
  <pubdate>1995</pubdate>
  <fpage>10</fpage>
  <lpage>-44</lpage>
</bibl>

<bibl id="B18">
  <title><p>Review of thermal properties of biological tissues</p></title>
  <aug>
    <au><snm>Giering</snm><fnm>K</fnm></au>
    <au><snm>Minet</snm><fnm>O</fnm></au>
    <au><snm>Lamprecht</snm><fnm>I</fnm></au>
    <au><snm>M{\"u}ller</snm><fnm>G</fnm></au>
  </aug>
  <source>Laser-induced interstitial therapy. SPIE Press, Bellingham,
  WA</source>
  <pubdate>1995</pubdate>
  <fpage>45</fpage>
  <lpage>-65</lpage>
</bibl>

<bibl id="B19">
  <title><p>Treatment planning for MRI-guided laser-induced interstitial
  thermotherapy of brain tumors—The role of blood perfusion</p></title>
  <aug>
    <au><snm>Schwarzmaier</snm><fnm>HJ</fnm></au>
    <au><snm>Yaroslavsky</snm><fnm>IV</fnm></au>
    <au><snm>Yaroslavsky</snm><fnm>AN</fnm></au>
    <au><snm>Fiedler</snm><fnm>V</fnm></au>
    <au><snm>Ulrich</snm><fnm>F</fnm></au>
    <au><snm>Kahn</snm><fnm>T</fnm></au>
  </aug>
  <source>Journal of Magnetic Resonance Imaging</source>
  <pubdate>1998</pubdate>
  <volume>8</volume>
  <issue>1</issue>
  <fpage>121</fpage>
  <lpage>127</lpage>
</bibl>

<bibl id="B20">
  <title><p>Identification of the Blood Perfusion Rate for Laser-Induced
  Thermotherapy in the Liver</p></title>
  <aug>
    <au><snm>Andres</snm><fnm>M</fnm></au>
    <au><snm>Blauth</snm><fnm>S</fnm></au>
    <au><snm>Leithäuser</snm><fnm>C</fnm></au>
    <au><snm>Siedow</snm><fnm>N</fnm></au>
  </aug>
  <pubdate>2019</pubdate>
</bibl>

</refgrp>
} 

\end{backmatter}
\end{document}